\documentclass[12pt]{article}
\usepackage{amssymb, amsmath, amsfonts,dsfont, mathrsfs,euscript,eufrak,colonequals}
\usepackage{textcomp}
\usepackage{latexsym}

\newtheorem{Theorem}{Theorem}
\newtheorem{Lemma}{Lemma}
\newtheorem{Proposition}{Proposition}
\newtheorem{Corollary}{Corollary}
\newtheorem{Remark}{Remark}
\newtheorem{Definition}{Definition}

\newcommand{\newparagraph}{\hspace{\parindent}}

\newcommand{\condA}{\mathrm{(A)}}
\newcommand{\condAone}{\mathrm{(A_1)}}
\newcommand{\condAtwo}{\mathrm{(A_2)}}
\newcommand{\condB}{\mathrm{(B)}}
\newcommand{\condC}{\mathrm{(C)}}

\newcommand{\N}{\mathbb N}
\newcommand{\R}{\mathbb R}

\newcommand{\Probab}{\mathds{P}}
\newcommand{\Expec}{ \mathds{E}\,}
\newcommand{\Disp}{\mathrm{\mathbf{Var}}\,}

\newcommand{\sign}{\mathrm{sign}}
\newcommand{\dd}{\,\mathrm{d}}
\newcommand{\id}{\mathds{1}}

\newcommand{\e}{\varepsilon}
\newcommand{\AppClass}{\mathcal{A}}
\newcommand{\APP}{\mathrm{APP}}

\newcommand{\SelfDecLaws}{\textit{\textbf{L}}}

\newcommand{\ContSet}{\mathbf{C}}
\newcommand{\lext}{\mathrm{lext}\,}
\newcommand{\rext}{\mathrm{rext}\,}
\newcommand{\mult}{\mathfrak{m}}
\newcommand{\StL}{S}
\newcommand{\EulerCon}{\mathscr{C}}

\newcommand{\nlambda}{\bar{\lambda}}

\newcommand{\svf}{\varphi}
\newcommand{\Laplsvf}{\widetilde{\svf}}
\newcommand{\deHaansvf}{\hat{\svf}}

\newcommand{\CorrFunc}{\mathcal{K}}
\newcommand{\EulerField}{\mathbb{E}}
\begin{document}
\author{A. A. Khartov}
\title{Asymptotic analysis of average case approximation complexity of Hilbert space valued random elements\footnote{The work was supported by
the Government of the Russian Federation megagrant 11.G34.31.0026, by JSC ``Gazprom Neft'', by the RFBR grant 13-01-00172, by the SPbGU grant 6.38.672.2013, and by the grant of Scientific school NSh-2504.2014.1.}}

\maketitle

\begin{abstract}
We study approximation properties of sequences of centered random elements $X_d$, $d\in\N$, with values in separable Hilbert spaces. We focus on sequences of tensor product-type and, in particular, degree-type random elements, which have covariance operators of corresponding tensor form. The average case approximation complexity $n^{X_d}(\e)$ is defined as the  minimal number of continuous linear functionals that is needed to approximate $X_d$ with relative $2$-average error not exceeding a given threshold $\e\in(0,1)$. In the paper we investigate $n^{X_d}(\e)$ for arbitrary fixed $\e\in(0,1)$ and $d\to\infty$.  Namely, we find criteria of (un)boundedness for $n^{X_d}(\e)$ on $d$ and of tending $n^{X_d}(\e)\to\infty$, $d\to\infty$, for any fixed $\e\in(0,1)$. In the latter case we obtain necessary and sufficient conditions for the following logarithmic asymptotics
\begin{eqnarray*}
\ln n^{X_d}(\e)= a_d+q(\e)b_d+o(b_d),\quad d\to\infty,
\end{eqnarray*}
at continuity points of a non-decreasing function $q\colon (0,1)\to\R$. Here $(a_d)_{d\in\N}$ is a sequence and $(b_d)_{d\in\N}$ is a positive sequence such that $b_d\to\infty$, $d\to\infty$.
Under rather weak assumptions, we show that for tensor product-type random elements only special quantiles of self-decomposable or, in particular, stable (for tensor degrees) probability distributions appear as functions $q$ in the asymptotics.

We apply our results to the tensor products of the Euler integrated processes with a given variation of smoothness parameters and to the tensor degrees of random elements with regularly varying eigenvalues of covariance operator.
\end{abstract}

\section{Introduction}\newparagraph
Let $X$ be a centered random element of some normed space $(Q, \|\cdot\|_Q)$. Let us approximate $X$ by the finite rank sums $\widetilde{X}^{(n)}=\sum_{k=1}^{n} l_k (X) \psi_k$, where $\psi_k$ are deterministic elements of $Q$ and $l_k$ are continuous linear functionals from the dual space $Q^*$. It is of theoretical and practical interest to make the relative average approximation error $(\Expec \|X-\widetilde{X}^{(n)}\|_Q^2/\Expec \|X\|_Q^2)^{1/2}$ smaller than a given error threshold $\e$ by choosing appropriate $n\in\N$ and optimal $\psi_k$ and $l_k$. Here we deal with \textit{linear approximation problem} in \textit{average case setting} (see \cite{Rit}). We call the minimal suitable value of such $n$, denoted by  $n^{X}(\e)$, the \textit{average case approximation complexity} of the  random element $X$ (see \cite{TraubWasWoz1}--\cite{TraubWers}).

Now we consider a sequence of centered random elements $X_d$, $d\in\N$, with values in normed spaces $(Q_d, \|\cdot\|_{Q_d})$, $d\in\N$, respectively. When $X_d$, $d\in\N$, are somehow related, it is of interest to look at the behaviour of the quantity $n^{X_d}(\e)$, as $d$ varies. In particular, for important \textit{linear tensor product approximation problems} (see \cite{NovWoz1}--\cite{NovWoz3}) we have  $Q_d=\otimes_{j=1}^d Q_{1,j}$ in appropriate sense, where $Q_{1,j}$, $j\in\N$, are some normed spaces. Here the covariances operators $K^{X_d}$ of  $X_d$, $d\in\N$, also have the appropriate tensor product form $K^{X_d}=\otimes_{j=1}^d K^{X_{1,j}}$, $d\in\N$, where $K^{X_{1,j}}$ is a covariance operator of given $Q_{1,j}$-valued centered random element $X_{1,j}$, $j\in\N$. Such  $X_d$ is called the \textit{tensor product of} $X_{1,1},\ldots, X_{1,d}$. If all $Q_{1,j}$ in  $Q_d$ are the same and all $K^{X_{1,j}}$ in $K^{X_d}$ are equal, then such $X_d$ is said to be the \textit{tensor degree of} $X_{1,1}$. A classical example of such objects is  the well known \textit{$d$-parametric Brownian sheet} (or  \textit{the Wiener--Chentsov random field}, see \cite{Adl} and \cite{Lifsh}). It, being considered as a random element of the space $L_2([0,1]^d)$, is a tensor degree of the Weiner process as a random element of $L_2([0,1])$. 

The described multivariate approximation problems find many applications in simulation methods, statistics, physics, and computational finance (see \cite{NovWoz0}). In this connection,  tractability questions becomes rather actual now. A sequence of approximation problems for $X_d, d\in\N$ is \textit{weakly tractable} if  $n^{X_d}(\e)$ is not exponential in $d$ or/and $\e^{-1}$. Otherwise, the sequence of the problems is \textit{intractable}. Special subclasses of weakly tractable problems are distinguished depending on the types of majorants for the quantity $n^{X_d}(\e)$ for all $d\in\N$ and $\e\in(0,1)$. There exist some results concerning the tractability of the described linear (tensor product) approximation problems, in which $Q_d$ are separable Hilbert spaces (see \cite{LifPapWoz1}, \cite{LifPapWoz2} and \cite{NovWoz1}--\cite{NovWoz3}). We will consider these problems within other less explored setting, namely, when $\e$ is arbitrarily close to zero but \textit{fixed} and $d$ goes to infinity. As noted in the book  \cite{NovWoz1} (see p. 6 and 289), such setting, being more appropriate for some models in computational finance, is also important. But the author is aware, in fact, of only one article \cite{LifTul} by M. A. Lifshits and E. V. Tulyakova on this subject. Our purpose is to complement their results.

The paper is organized as follows. In Section 2 we provide a formal problem setting. In Section 3 we consider sequences of random elements $X_d$, $d\in\N$, of separable Hilbert spaces. In particular, we find necessary and sufficient conditions that for almost all $\e\in(0,1)$ the quantity $n^{X_d}(\e)$ has a special form of logarithmic asymptotics as $d\to\infty$. In Sections 4 and 5 the same problem is solved for the sequences of tensor product-type and degree-type random elements, respectively. We show that, under rather weak assumption, $n^{X_d}(\e)$ depends on $\e$ according to some \textit{self-decomposable} or, in particular, \textit{stable probability distribution}. We apply the obtained criteria to tensor products of the Euler integrated processes with a given variation of the smothness parameters and to tensor degrees of random elements with a given regular variation of eigenvalues of covariance operator. In Appendix we provide necessary facts from probability theory about self-decomposable and stable distributions and also about related limit theorems, which are the main tools of our work.

Throughout the article, we use the following notation. We write $a_n\sim b_n$ iff $a_n/b_n\to 1$, $n\to\infty$. We denote by $\N$ and $\R$ the sets of positive integer and real numbers, respectively. We set $\ln_+ x\colonequals \max\{1, \ln x\}$ for all $x>0$. The quantity $\id(A)$ equals one for the true logic propositions $A$ and zero for the false ones. We always use $\|\cdot\|_B$ for the norm, which some space $B$ is equipped with. For any function $f$ we will denote by $\ContSet(f)$ the set of all its continuity points and by $f^{-1}$ the generalized inverse function $f^{-1}(y)\colonequals\inf\bigl\{x\in\R: f(x)\geqslant y\bigr\}$, where $y$ is from the range of $f$. By \textit{distribution function} $F$ we mean the non-decreasing function $F$ on $\R$ that is right-continuous on $\R$, $\lim\limits_{x\to-\infty} F(x)=0$, and $\lim\limits_{x\to\infty} F(x)=1$. Following \cite{LinOstr}, the boundaries of growth points of distribution function $F$ will be denoted by
\begin{eqnarray*}
\lext F\colonequals\inf\{x\in\R: F(x)>0\}\geqslant -\infty,\qquad\rext F\colonequals\sup\{x\in\R: F(x)<1\}\leqslant \infty.
\end{eqnarray*}
A distribution function $F$ is called \textit{degenerate} if   $F(x)=\id(x\geqslant a)$ for any $x$ and some constant $a$.

\section{Basic definitions and problem setting}\newparagraph
We consider a sequence of random elements $X_d$, $d\in\N$, with values in separable Hilbert spaces $H_d$, $d\in\N$, respectively. Assume that every $X_d$ has zero mean and $\Expec \|X_d\|_{H_d}^2<\infty$, $d\in\N$. We will investigate the \textit{average case approximation complexity} (simply the \textit{approximation complexity} for short) of $X_d$, $d\in\N$:
\begin{eqnarray}\label{def_nXde}
n^{X_d}(\e)\colonequals\min\bigl\{n\in\N:\, e^{X_d}(n)\leqslant \e\, e^{X_d}(0)  \bigr \},
\end{eqnarray}
where $\e\in(0,1)$ is a given error threshold, and
\begin{eqnarray*}
e^{X_d}(n)\colonequals\inf\Bigl\{ \bigl(\Expec\bigl\|X_d - \widetilde X^{( n)}_d\bigr\|_{H_d}^2\bigr)^{1/2} :  \widetilde X^{(n)}_d\in \AppClass_n^{X_d}\Bigr\}
\end{eqnarray*}
is the smallest 2-average error among all linear approximations of $X_d$, $d\in\N$, having rank $n\in\N$. The corresponding classes  of linear algorithms are denoted by
\begin{eqnarray*}
\AppClass_n^{X_d}\colonequals \Bigl\{\sum_{m=1}^{n}l_m(X_d)\,\psi_m :  \psi_m \in H_d,\, l_m\in H_d^*  \Bigr \}, \quad d\in\N,\quad n\in\N.
\end{eqnarray*}
We work with \textit{relative errors}, thus taking into account the following ``size'' of $X_d$:
\begin{eqnarray*}
e^{X_d}(0)\colonequals  \bigl(\Expec\|X_d\|_{H_d}^2\bigr)^{1/2},
\end{eqnarray*}
which is the approximation error of $X_d$ by zero element of $H_d$. The approximation complexity $n^{X_d}(\e)$ is considered as a function depending on two variables $d$ and $\e$. The general goal is to understand the character of this dependence for the given sequence $(X_d)_{d\in\N}$. 

The linear tensor product approximation problems are of our particular interest. We will study these only within the following construction. We suppose that every $H_d$ is defined by the Hilbertian tensor product $H_d\colonequals \otimes_{j=1}^d H_{1,j}$, where every $H_{1,j}$ is a given separable Hilbert space, $j\in\N$. We also suppose that for every $d$ the random element $X_d$ is a \textit{tensor product} of some zero-mean $H_{1,j}$-valued random elements $X_{1,j}$, $j=1,\ldots,d$, i.e. $X_d$ has zero mean and the covariance operator $K^{X_d}\colonequals\otimes_{j=1}^d K^{X_{1,j}}$, where $K^{X_{1,j}}$ is a covariance operator of $X_{1,j}$ for every $j\in\N$. Following \cite{KarNazNik} and \cite{LifZani}, for such $X_d$ we will write $X_d=\otimes_{j=1}^d X_{1,j}$ for short. In particular, if $H_{1,j}$ are the same and $K^{X_{1,1}}=\ldots=K^{X_{1,d}}$, then the element $X_d$ is called the \textit{tensor degree} of $X_{1,1}$ and we will write $X_d=X_{1,1}^{\otimes d}$ in this case.

We now turn to tractability of approximation problems. Let us give basic definitions related to the tractability in context of our setting \eqref{def_nXde}. Let $\APP$ denote a sequence of linear approximation problems for $X_d$, $d\in\N$. We say that
\begin{itemize}
\item $\APP$ is \textit{weakly tractable} iff
\begin{eqnarray}\label{def_WT}
\lim_{d+\e^{-1}\to\infty}\dfrac{\ln_+ n^{X_d}(\e)}{d+\e^{-1}}=0;
\end{eqnarray}
\item $\APP$ is \textit{quasi-polynomially tractable} iff there are numbers $C>0$ and $s\geqslant0$ such that
\begin{eqnarray}\label{def_QPT}
n^{X_d}(\e)\leqslant C\,\exp\bigl\{s(1+\ln \e^{-1})\ln_+ d\bigr\}\quad\text{for all}\quad d\in\N, \,\,\e\in(0,1);
\end{eqnarray}
\item $\APP$ is \textit{polynomially tractable} iff there are numbers  $C>0$, $s\geqslant 0$, and $p\geqslant0$ such that 
\begin{eqnarray}\label{def_PT}
n^{X_d}(\e)\leqslant C\,\e^{-s}\,d^{\,p} \quad\text{for all}\quad d\in\N, \,\,\e\in(0,1);
\end{eqnarray}
\item $\APP$ is \textit{strong polynomially tractable} iff there are numbers $C>0$ and $s\geqslant0$ such that
\begin{eqnarray}\label{def_SPT}
n^{X_d}(\e)\leqslant C\,\e^{-s}\quad\text{for all}\quad d\in\N, \,\,\e\in(0,1).
\end{eqnarray}
\end{itemize}
If $\APP$ is not weakly tractable, then it is called \textit{intractable}.  In particular, if $n^{X_d}(\e)$ increases at least exponentially in $d$, then we say that $\APP$ has a \textit{curse of dimensionality}.

Let us give a short review of some results concerning the tractability in our setting. For the  linear approximation problems there exist criteria of each type of tractability in terms of eigenvalues of the operators $K^{X_d}$, $d\in\N$. These results can be found in \cite{NovWoz1} (see p. 245, 256). The described linear tensor product approximation problems were investigated in the recent paper \cite{LifPapWoz1}, where the necessary and sufficient conditions for all types of tractability were given in terms of eigenvalues of the  covariance operators $K^{X_{1,j}}$, $j\in\N$. Also the tensor product-type Korobov kernels were studied there. In \cite{LifPapWoz2} the authors obtained  results concerning tractability of tensor products of the Euler and Wiener integrated processes with varying smoothness parameters. 

In the paper we investigate the approximation complexity $n^{X_d}(\e)$ of general and tensor product-type random elements for arbitrary fixed $\e\in(0,1)$ and $d\to\infty$. Namely, we search  necessary and sufficient conditions that $n^{X_d}(\e)$ has the following logarithmic asymptotics 
\begin{eqnarray}\label{conc_nXde_LogAsymp_form}
\ln n^{X_d}(\e)= a_d+q(\e)b_d+o(b_d),\quad d\to\infty,\quad\text{for all}\quad \e\in\ContSet(q),
\end{eqnarray}
where $(a_d)_{d\in\N}$ is a sequence, $(b_d)_{d\in\N}$ is a positive sequence such that $b_d\to\infty$, $d\to\infty$, and the function $q\colon (0,1)\to\R$ is non-increasing. To consider the asymptotics of the form \eqref{conc_nXde_LogAsymp_form} is quite natural and, moreover, as we will see below, it is inherent in the linear tensor product approximation problems. The author is not aware of any results, which could provide the full solution in such setting. But there exist three papers \cite{LifTul}--\cite{LifZani2014} concerning asymptotic analysis of $n^{X_d}(\e)$ for fixed $\e$ and $d\to\infty$, where only the first one deals with objects of our interest. In \cite{LifTul} the authors considered only the tensor degree-type random elements under some assumptions. We will discuss the corresponding theorem in Section 5.

As in case of tractability, eigenvalues of covariance operators of $X_d$, $d\in\N$ or $X_{1,j}$, $j\in\N$, play a crucial role in asymptotic analisys of the quantity $n^{X_d}(\e)$ as $d\to\infty$. For convenience, throughout the paper we will use the following unified notation for the correlation characteristics. For a given Hilbertian random element $Z$ we will denote by $K^Z$ its covariance operator. The sequences $(\lambda_k^Z)_{k\in\N}$ and $(\psi^Z_k)_{k\in\N}$ will denote the non-increasing sequence of eigenvalues and the corresponding  sequence of eigenvectors of $K^Z$, respectively, i.e. $K^Z \psi^Z_k=\lambda^Z_k\psi^Z_k$, $k\in\N$. If $Z$ is a random element of $p$-dimensional space, then we formally set $\lambda_k^Z\colonequals0$, and $\psi_k^Z\colonequals0$ for $k>p$. The trace of $K^Z$ will be denoted by $\Lambda^{Z}$, thus $\Lambda^{Z}=\sum_{k=1}^\infty \lambda^{Z}_k$. We will also use the notation $\nlambda^Z_k\colonequals \lambda^Z_k/\Lambda^Z$, $k\in\N$, (i.e. $\sum_{k=1}^\infty \nlambda^{Z}_k=1$) and  $\mult^Z(x)\colonequals\sum_{k=1}^\infty \id\bigl(\nlambda^Z_k=x\bigr)$ for any $x>0$.

\section{Approximation of general random elements}\newparagraph
Here we consider a general sequence of random elements $X_d$, $d\in\N$, of abstract separable Hilbert spaces $H_d$, $d\in\N$, respectively, without any assumptions on the spectral structure of the corresponding covariance operators $K^{X_d}$, $d\in\N$. Let $H_d$ be equipped with inner product $(\,\cdot\,,\,\cdot\,)_{H_d}$.  We always assume that every $X_d$ has zero mean and satisfies $\Expec \|X_d\|_{H_d}^2<\infty$, $d\in\N$. Then self-adjoint non-negative definite operators $K^{X_d}$, $d\in\N$, have the finite traces:
\begin{eqnarray}\label{conc_LambdaXd_GRE}
\Lambda^{X_d}=\sum_{k=1}^\infty \lambda^{X_d}_k=\Expec \|X_d\|_{H_d}^2=e^{X_d}(0)^2<\infty,\quad d\in\N.
\end{eqnarray}
To omit the pathological cases, we always assume that $\lambda_1^{X_d}>0$ for all $d\in\N$.

It is well known (see \cite{WasWoz}) that for any $n\in\N$ the following random element
\begin{eqnarray}\label{def_Xdn_GRE}
\widetilde X^{(n)}_d\colonequals\sum_{k=1}^n (X_d,\psi^{X_d}_k)_{H_d}\, \psi^{X_d}_k\in\AppClass_n^{X_d}
\end{eqnarray}
minimizes the 2-average case error. Hence formula \eqref{def_nXde} is reduced to
\begin{eqnarray*}
n^{X_d}(\e)=\min\Bigl\{n\in\N:\, \Expec 
\bigl\|X_d-\widetilde{X}^{(n)}_d\bigr\|_{H_d}^2\leqslant\e^2\,\Expec \|X_d\|_{H_d}^2 \Bigr\}.
\end{eqnarray*}
From \eqref{conc_LambdaXd_GRE} and \eqref{def_Xdn_GRE} we infer the following representation of the approximation complexity
\begin{eqnarray*}
n^{X_d}(\e)=\min\Bigl\{n\in\N:\, \sum_{k=n+1}^\infty \lambda^{X_d}_k\leqslant\e^2\,\Lambda^{X_d} \Bigr\}.
\end{eqnarray*}
In terms of $\nlambda^{X_d}_k$, $k\in\N$, it takes the form 
\begin{eqnarray}
n^{X_d}(\e)&=&\min\Bigl\{n\in\N:\, \sum_{k=n+1}^\infty \nlambda^{X_d}_k\leqslant\e^2 \Bigr\}\label{conc_nXde1}\\
&=&\min\Bigl\{n\in\N: \sum_{k=1}^n \nlambda^{X_d}_k\geqslant 1-\e^2\Bigr\}\label{conc_nXde2}.
\end{eqnarray}

\subsection{Boundedness of the approximation complexity}\newparagraph
Before proceeding to the asymptotic analysis of the quantity $n^{X_d}(\e)$, we find criteria of its boundedness and unboundedness on $d$ for any fixed $\e\in(0,1)$. The next simple proposition provides the conditions for $n^{X_d}(\e)\to\infty$ as $d\to\infty$ in terms of the first normed eigenvalues $\nlambda^{X_d}_1$, $d\in\N$.
\begin{Proposition}\label{pr_nXdetoinfty_GRE}
The following conditions are equivalent:
\begin{itemize}
\item[$(i)$]\quad $\lim\limits_{d\to\infty}n^{X_d}(\e)=\infty$\quad for all\quad $\e\in(0,1);$
\item[$(ii)$]\quad $\lim\limits_{d\to\infty}\nlambda^{X_d}_1=0.$
\end{itemize}
\end{Proposition}
\textbf{Proof of Proposition \ref{pr_nXdetoinfty_GRE}.}\quad $(i)\Rightarrow(ii)$. Suppose that, contrary to our claim, there exists a subsequence $(d_l)_{l\in\N}$ such that $c\colonequals\lim\limits_{l\to\infty} \nlambda^{X_{d_l}}_1>0$. Choose $\e\in(0,1)$ to satisfy $1-\e^2<c$. Then, according to \eqref{conc_nXde2}, we have $n^{X_{d_l}}(\e)=1$ for all sufficiently large $l\in\N$. This contradicts our assumption $(i)$.

$(i)\Leftarrow(ii)$. For all $\e\in(0,1)$ the statement $n^{X_d}(\e)\to\infty$, $d\to\infty$,
follows from the formula \eqref{conc_nXde2} and the inequality:
\begin{eqnarray}\label{conc_nXde_lambdaXd1_ineq_GRE}
n^{X_d}(\e)\geqslant\sum_{k=1}^{n^{X_d}(\e)}\dfrac{\nlambda^{X_d}_k}{\nlambda^{X_d}_1}\geqslant \dfrac{1-\e^2}{\nlambda^{X_d}_1}.\quad \Box
\end{eqnarray}

As a consequence of the previous proposition, we obtain the following criterion of unboundedness of $n^{X_d}(\e)$ on $d\in\N$ for any fixed $\e\in(0,1)$.

\begin{Proposition}\label{pr_nXdeunbounded_GRE}
The following conditions are equivalent:
\begin{itemize}
\item[$(i)$]\quad $\sup\limits_{d\in\N} n^{X_d}(\e)=\infty$\quad for all\quad $\e\in(0,1);$ 
\item[$(ii)$]\quad $\inf\limits_{d\in\N}\nlambda^{X_d}_1=0.$
\end{itemize}
\end{Proposition}

Let us formulate the criterion of boundedness of the approximation complexity $n^{X_d}(\e)$ on $d\in\N$ for any fixed $\e\in(0,1)$.

\begin{Proposition}\label{pr_nXdebounded_GRE}
The following conditions are equivalent:
\begin{itemize}
\item[$(i)$] \quad $\sup\limits_{d\in\N}n^{X_d}(\e)<\infty$\quad for all\quad $\e\in(0,1);$
\item[$(ii)$]\quad $\lim\limits_{t\to 0}\,\sup\limits_{d\in\N}\sum\limits_{k=1}^\infty\nlambda^{X_d}_k\id\bigl(\nlambda^{X_d}_k<t\bigr) =0.$
\end{itemize}
\end{Proposition}

Before proving this proposition, we introduce the following auxiliary quantity
\begin{eqnarray}\label{def_nlambdaXde}
\nlambda^{X_d}(\e)\colonequals\nlambda^{X_d}_{n^{X_d}(\e)},\quad d\in\N,\quad \e\in(0,1),
\end{eqnarray}
which admits the following representation
\begin{eqnarray}\label{conc_nlambdaXde_sup}
	\nlambda^{X_d}(\e)
	=\sup\Bigl\{x\in\R: \sum_{k=1}^\infty\nlambda^{X_d}_k\,\id\bigr(\nlambda^{X_d}_k< x\bigl)\leqslant \e^2\Bigr\}.
\end{eqnarray}
Obtain useful inequalities connecting $\nlambda^{X_d}(\e)$ and $n^{X_d}(\e)$. First, from the definitions of these quantities we conclude that
\begin{eqnarray}\label{conc_nXde_nlambdaXde_upbound}
n^{X_d}(\e)\leqslant \nlambda^{X_d}(\e)^{-1}\quad\text{for all}\quad d\in\N, \quad\e\in(0,1).
\end{eqnarray}
Next, for any  $\e_1\in(0,1)$, $\e_2\in(\e_1,1)$, and $d\in\N$ we have
\begin{eqnarray*}
	n^{X_d}(\e_1)\geqslant n^{X_d}(\e_1)-n^{X_d}(\e_2)+1\geqslant \sum_{k=n^{X_d}(\e_2)+1}^{n^{X_d}(\e_1)} \dfrac{\nlambda^{X_d}_k}{\nlambda^{X_d}(\e_2)}+1,
\end{eqnarray*} 
where
\begin{eqnarray*}
	\sum_{k=n^{X_d}(\e_2)+1}^{n^{X_d}(\e_1)} \nlambda^{X_d}_k
	&=&\sum_{k=1}^{n^{X_d}(\e_1)} \nlambda^{X_d}_k-\sum_{k=1}^{n^{X_d}(\e_2)} \nlambda^{X_d}_k\\
	&\geqslant&(1-\e_1^2)- \bigl((1-\e_2^2)+\nlambda^{X_d}(\e_2)\bigr)\\
	&=&\e_2^2-\e_1^2- \nlambda^{X_d}(\e_2).
\end{eqnarray*}
Finally, we get
\begin{eqnarray}\label{conc_nXde_nlambdaXde_lowbound}
n^{X_d}(\e_1)\geqslant \dfrac{\e_2^2-\e_1^2}{\nlambda^{X_d}(\e_2)}\quad\text{for all}\quad d\in\N, \,\,\e_1\in(0,1), \,\,\e_2\in(\e_1,1).
\end{eqnarray}

\noindent\textbf{Proof of Proposition \ref{pr_nXdebounded_GRE}.}\quad  First, we show that the condition $(i)$ is equivalent to
\begin{eqnarray}\label{conc_nlambdaXde_zero}
\inf\limits_{d\in\N}\nlambda^{X_d}(\e)>0\quad \text{for all}\quad\e\in(0,1).
\end{eqnarray}
Indeed, from \eqref{conc_nXde_nlambdaXde_upbound} we conclude the implication $\eqref{conc_nlambdaXde_zero}\Rightarrow(i)$. Next, let us fix any $\e\in(0,1)$ and $c\in(0,1)$. Using inequality \eqref{conc_nXde_nlambdaXde_lowbound} with $\e_1=c\e$ and $\e_2=\e$, we obtain
\begin{eqnarray*}
\nlambda^{X_d}(\e)\geqslant \dfrac{(1-c^2)\e^2}{\sup\limits_{d\in\N} n^{X_d}(c\e)}\quad\text{for all}\quad d\in\N.
\end{eqnarray*}
From this we get the implication $(i)\Rightarrow\eqref{conc_nlambdaXde_zero}$.

It only remains to verify $(ii)\Leftrightarrow\eqref{conc_nlambdaXde_zero}$. But it follows from the assertion that for any $\e\in(0,1)$ and some $t_\e\in(0,1)$ we have $\inf_{d\in\N}\nlambda^{X_d}(\e)\geqslant t_\e$ iff $\sup_{d\in\N}\sum_{k=1}^\infty\nlambda^{X_d}_k\,\id\bigl(\nlambda^{X_d}_k<t_\e\bigr)\leqslant \e^2$. \quad $\Box$

\subsection{Logarithmic asymptotics of the approximation complexity}\newparagraph
Here we start asymptotic analysis of the approximation complexity $n^{X_d}(\e)$. Our next theorem establishes a connection between the asymptotics of the form \eqref{conc_nXde_LogAsymp_form}, for given $(a_d)_{d\in\N}$ and $(b_d)_{d\in\N}$, and the convergence of the following distribution functions
\begin{eqnarray}\label{def_GXdadbd_GRE}
G^{X_d}_{a_d,b_d}(x)\colonequals \sum_{k=1}^{\infty} \nlambda^{X_d}_k\, \id\bigl( \nlambda^{X_d}_k\geqslant e^{-a_d-xb_d}\bigr),\quad x\in\R,\quad d\in\N.
\end{eqnarray}

\begin{Theorem}\label{th_nXde_LogAsymp_GRE}
Let $(a_d)_{d\in\N}$ be a sequence, $(b_d)_{d\in\N}$ be a positive sequence such that $b_d\to\infty$, $d\to\infty$. Let a non-increasing function $q\colon (0,1)\to\R$ and a distribution function $G$ satisfy the equation   $q(\e)=G^{-1}(1-\e^2)$ for all $\e\in\ContSet(q)$. For the asymptotics
\begin{eqnarray}\label{th_nXde_LogAsymp_GRE_nXde_asymp}
\ln n^{X_d}(\e)=a_d+q(\e)b_d+o(b_d),\quad d\to\infty,\quad\text{for all}\quad \e\in\ContSet(q),
\end{eqnarray}
it is necessary and sufficient to have the weak convergence\footnote{We follow the definition from \cite{Petr}, p. 16. Equivalent definitions can be found in \cite{Fell}, p. 248--251.} $G^{X_d}_{a_d,b_d}\Rightarrow G$ as $d\to\infty$, i.e.
\begin{eqnarray}\label{th_nXde_LogAsymp_GRE_GXdadbd_conv}
\lim_{d\to\infty} G^{X_d}_{a_d,b_d}(x)= G(x) \quad\text{for all}\quad x\in \ContSet(G).
\end{eqnarray}
\end{Theorem}
\textbf{Proof of Theorem \ref{th_nXde_LogAsymp_GRE}.}\quad It is well known that condition \eqref{th_nXde_LogAsymp_GRE_GXdadbd_conv} is equivalent to the convergence $\lim\limits_{d\to\infty}\bigl(G^{X_d}_{a_d,b_d}\bigr)^{-1}(y)= G^{-1}(y)$ for all continuity points $y$ of $G^{-1}$ (see \cite{Vaart}, p. 305). By theorem of continuity of composite functions, we have $\e \in \ContSet(q) \Leftrightarrow 1-\e^2\in\ContSet(G^{-1})$. Therefore \eqref{th_nXde_LogAsymp_GRE_GXdadbd_conv} is equivalent to the following condition
\begin{eqnarray*}
\lim\limits_{d\to\infty}\bigl(G^{X_d}_{a_d,b_d}\bigr)^{-1}(1-\e^2)= q(\e)\quad\text{for all}\quad \e\in\ContSet(q).
\end{eqnarray*}
According to \eqref{def_nlambdaXde} and \eqref{conc_nlambdaXde_sup}, we note that
\begin{eqnarray}
\bigl|\ln\nlambda^{X_d}(\e)\bigr|
&=& \inf\Bigl\{y\in\R: \sum_{k=1}^\infty\nlambda^{X_d}_k\,\id\bigr(\nlambda^{X_d}_k\geqslant e^{-y}\bigl)\geqslant 1-\e^2\Bigr\}\label{conc_nlambdaXde_inf}\\
&=& \inf\Bigl\{y\in\R: G^{X_d}_{a_d,b_d}\bigl((y-a_d)/b_d\bigr)\geqslant 1-\e^2\Bigr\}\nonumber\\
&=& a_d+b_d\cdot\inf\bigl\{z\in\R: G^{X_d}_{a_d,b_d}(z)\geqslant 1-\e^2\bigr\}\nonumber\\
&=& a_d+ b_d\cdot \bigl(G^{X_d}_{a_d,b_d}\bigr)^{-1}(1-\e^2).\nonumber
\end{eqnarray}
Thus \eqref{th_nXde_LogAsymp_GRE_GXdadbd_conv} is equivalent to the following condition  
\begin{eqnarray}\label{conc_nlambdaXde_LogAsymp_GRE}
\bigl|\ln\nlambda^{X_d}(\e)\bigr|=a_d+q(\e)b_d+o(b_d),\quad d\to\infty,\quad\text{for all}\quad \e\in\ContSet(q).
\end{eqnarray}

Let us prove the implication \eqref{conc_nlambdaXde_LogAsymp_GRE} $\Rightarrow$ \eqref{th_nXde_LogAsymp_GRE_nXde_asymp} . Fix $\e\in\ContSet(q)$ and arbitrarily small $h>0$. From  \eqref{conc_nlambdaXde_LogAsymp_GRE} and \eqref{conc_nXde_nlambdaXde_upbound}  we have the inequality $\ln n^{X_d}(\e)\leqslant a_d+q(\e)b_d+h b_d$ for all large enough $d\in\N$. Since the function $q$ is non-increasing, the set $\ContSet(q)$ is dense in the interval $(0,1)$. There exists $\tau_1\in(1,1/\e)$ such that $\tau_1\e\in\ContSet(q)$ and $q(\tau_1\e)-q(\e)\geqslant-h$. By the inequality \eqref{conc_nXde_nlambdaXde_lowbound}, we have 
\begin{eqnarray*}
\ln n^{X_d}(\e)\geqslant\bigl|\ln\nlambda^{X_d}(\tau_1\e)\bigr|+ \ln \bigl((\tau_1\e)^2-\e^2\bigr).
\end{eqnarray*}
According to \eqref{conc_nlambdaXde_LogAsymp_GRE}, for all large enough $d\in\N$ we obtain
\begin{eqnarray*}
\ln n^{X_d}(\e)&\geqslant& a_d+q(\tau_1\e)b_d-h b_d+\ln \bigl((\tau_1\e)^2-\e^2\bigr)\\
&\geqslant& a_d+q(\e)b_d-2hb_d+\ln \bigl((\tau_1\e)^2-\e^2\bigr).
\end{eqnarray*}
Since $b_d\to\infty$ as $d\to\infty$, it follows that $\ln n^{X_d}(\e)\geqslant a_d+q(\e)b_d-3hb_d$ for all large enough $d\in\N$. Thus the asymptotics \eqref{th_nXde_LogAsymp_GRE_nXde_asymp} follows from the obtained estimates for $n^{X_d}(\e)$.

Prove \eqref{th_nXde_LogAsymp_GRE_nXde_asymp} $\Rightarrow$ \eqref{conc_nlambdaXde_LogAsymp_GRE}. Fix $\e\in\ContSet(q)$ and arbitrarily small $h>0$. From  \eqref{th_nXde_LogAsymp_GRE_nXde_asymp} and \eqref{conc_nXde_nlambdaXde_upbound} we obtain the inequality $\bigl|\ln\nlambda^{X_d}(\e)\bigr|\geqslant a_d+q(\e)b_d-h b_d$ for all large enough $d\in\N$. Find $\tau_2\in(0,1)$ such that $\tau_2\e\in\ContSet(q)$ and $q(\tau_2\e)-q(\e)\leqslant h$. The inequality \eqref{conc_nXde_nlambdaXde_lowbound} yields 
\begin{eqnarray*}
\bigl|\ln\nlambda^{X_d}(\e)\bigr|\leqslant\ln n^{X_d}(\tau_2\e)-\ln \bigl(\e^2-(\tau_2\e)^2\bigr).
\end{eqnarray*}
On account of the asymptotics \eqref{th_nXde_LogAsymp_GRE_nXde_asymp}, for all large enough $d\in\N$ we obtain
\begin{eqnarray*}
\bigl|\ln\nlambda^{X_d}(\e)\bigr|&\leqslant& a_d+q(\tau_2\e)b_d+h b_d-\ln \bigl(\e^2-(\tau_2\e)^2\bigr)\\
&\leqslant& a_d+q(\e)b_d+2hb_d-\ln \bigl(\e^2-(\tau_2\e)^2\bigr).
\end{eqnarray*}
Since $b_d\to\infty$ as $d\to\infty$, we have $\bigl|\ln\nlambda^{X_d}(\e)\bigr|\leqslant a_d+q(\e)b_d+3hb_d$ for all large enough $d\in\N$. The obtained estimates yield the required asymptotics \eqref{conc_nlambdaXde_LogAsymp_GRE}.\quad $\Box$\\

Let us mention some elementary facts for the functions $G$ and $q$ from Theorem \ref{th_nXde_LogAsymp_GRE}.

\begin{Remark}
If for any given $\e\in(0,1)$ the distribution function $G$ strongly increases on the right $(\text{left}\,)$ neighbourhood of $q(\e)$, then $q$ is left-$(\text{right-})$continuous at $\e$.
\end{Remark}

Indeed, suppose that $G$ strongly increases on the right neighbourhood of $q(\e)$ (the left case is similar). Then for any $\delta>0$ we have $G(q(\e)+\delta)>1-\e^2$
and, consequently, there exists $\tau_\delta>0$ such that $G(q(\e)+\delta)>1-(\e-\tau)^2$ for all $\tau\in(0,\tau_\delta)$. Applying $G^{-1}$ to the previous inequality, we obtain $q(\e-\tau)-q(\e)\leqslant\delta$, which gives left-continuity of $q$ at $\e$ by monotony of $q$.

As a consequence of this remark, we provide the following useful note.
\begin{Remark}\label{rem_Contq_GRE}
If the distribution function $G$ is degenerate or strongly increases on the non-empty interval $(\lext G, \rext G)$, then $q$ is continuous on $(0,1)$.
\end{Remark}

\section{Approximation of tensor product-type random elements}\newparagraph
In this section we consider sequences of tensor product-type random elements. Suppose that we have a sequence of zero-mean random elements $X_{1,j}$, $j\in\N$, of separable Hilbert spaces $H_{1,j}$, $j\in\N$, respectively. Assume that every $X_{1,j}$ satisfies $\Expec \|X_{1,j}\|^2_{H_{1,j}}<\infty$, $j\in\N$. Let $X_d=\otimes_{j=1}^d X_{1,j}$, $d\in\N$. It is well known that eigenvalues and eigenvectors of covariance operator $K^{X_d}$ have the following multiplicative form:
\begin{eqnarray}\label{conc_eigenpairs_TPRE}
\prod_{j=1}^d\lambda^{X_{1,j}}_{k_j},\quad \bigotimes_{j=1}^d\psi^{X_{1,j}}_{k_j},\qquad k_1, k_2, \ldots,k_d\in\N.
\end{eqnarray}
Hence for traces $\Lambda^{X_d}$  of $K^{X_d}$, $d\in\N$, we have the formula
\begin{eqnarray}\label{conc_LambdaXd_TPRE}
\Lambda^{X_d}=\sum_{k_1, k_2, \ldots,k_d\in\N} \prod_{j=1}^d\lambda^{X_{1,j}}_{k_j}=\prod_{j=1}^d\sum_{i=1}^\infty\lambda^{X_{1,j}}_i=\prod_{j=1}^d\Lambda^{X_{1,j}}.
\end{eqnarray}
Throughout this section we assume that $\lambda^{X_{1,j}}_1>0$ for all $j\in\N$. Since the first (the maximal) eigenvalue $\lambda^{X_d}_1$ of $K^{X_d}$ equals $\prod_{j=1}^d\lambda^{X_{1,j}}_1$, this convention is equivalent to the assumption from Section~3, namely $\lambda^{X_d}_1>0$, $d\in\N$, .

\subsection{Boundedness of the approximation complexity}\newparagraph
By analogy with Section 3, we first consider the boundedness conditions of the approximation complexity $n^{X_d}(\e)$ on $d$ for any fixed $\e\in(0,1)$. For described tensor product-type random elements $X_d$, $d\in\N$, the following propositions show that for any fixed $\e\in(0,1)$ either the quantity $n^{X_d}(\e)$ tends to infinity as $d\to\infty$ or it is a bounded function on $d\in\N$. 

\begin{Proposition}\label{pr_nXdetoinfty_TPRE}
The following conditions are equivalent:
\begin{itemize}
\item[$(i)$]\quad  $\lim\limits_{d\to\infty}n^{X_d}(\e)=\infty$\quad for all\quad $\e\in(0,1);$
\item[$(ii)$]\quad $\sum\limits_{j=1}^{\infty} \sum\limits_{k=2}^\infty\dfrac{\lambda^{X_{1,j}}_k}{\lambda^{X_{1,j}}_1}= \infty.$
\end{itemize}
\end{Proposition}
\textbf{Proof of Proposition \ref{pr_nXdetoinfty_TPRE}.}\quad By Proposition \ref{pr_nXdetoinfty_GRE}, the relation $(i)$ is equivalent to the convergence $\lim\limits_{d\to\infty}(\Lambda^{X_d}/\lambda^{X_d}_1)=\infty$, where, as it is easily seen:
\begin{eqnarray}\label{conc_lambdaXd1_TPRE}
\dfrac{\Lambda^{X_d}}{\lambda^{X_d}_1}=\prod_{j=1}^d \dfrac{\Lambda^{X_{1,j}}}{\lambda^{X_{1,j}}_1}=\prod_{j=1}^{d} \Biggl(1+ \sum_{k=2}^\infty\dfrac{\lambda^{X_{1,j}}_k}{\lambda^{X_{1,j}}_1}\Biggr).
\end{eqnarray}
The last product goes to infinity as $d\to\infty$ iff the relation $(ii)$ holds.\quad $\Box$

\begin{Proposition}\label{pr_nXdebounded_TPRE}
The following conditions are equivalent:
\begin{itemize}
\item[$(i)$] \quad $\sup\limits_{d\in\N} n^{X_d}(\e)<\infty$\quad for all \quad $\e\in(0,1);$
\item[$(ii)$] \quad $\sum\limits_{j=1}^{\infty} \sum\limits_{k=2}^\infty\dfrac{\lambda^{X_{1,j}}_k}{\lambda^{X_{1,j}}_1}<\infty.$
\end{itemize}
\end{Proposition}
\textbf{Proof of Proposition \ref{pr_nXdebounded_TPRE}.}\quad First, we note that $(ii)$ is equivalent to 
\begin{eqnarray*}
C\colonequals\sup_{d\in\N} \bigl(\Lambda^{X_d}/\lambda^{X_d}_1\bigr)<\infty.
\end{eqnarray*}
Under the condition $(i)$, the last directly follows from Proposition \ref{pr_nXdeunbounded_GRE}. It only remains to check that the assumption $C<\infty$ implies $(i)$. Let $\hat\lambda^{X_\infty}_k$, $k\in\N$ be renumbered sequence of numbers $\prod_{j=1}^\infty \bigl(\lambda^{X_{1,j}}_{k_j}/\lambda^{X_{1,j}}_1\bigr)$, where $k_j\in\N$ and $\lim\limits_{j\to\infty}k_j=1$. According to the multiplicative structure of $(\lambda^{X_d}_k)_{k\in\N}$ (see \eqref{conc_eigenpairs_TPRE}), we have
\begin{eqnarray*}
\sup_{d\in\N}\sum\limits_{k=1}^\infty \nlambda^{X_d}_k\,\id\bigl(\nlambda^{X_d}_k<t \bigr)
&\leqslant& \sup_{d\in\N}\sum\limits_{k=1}^\infty \bigl(\lambda^{X_d}_k/\lambda^{X_d}_1\bigr)\,\id\bigl(\lambda^{X_d}_k/\lambda^{X_d}_1<C t \bigr)\\
&=& \sup_{d\in\N}\sum\limits_{k\in\N^d} \prod_{j=1}^d \bigl(\lambda^{X_{1,j}}_{k_j}/\lambda^{X_{1,j}}_1\bigr)\,\id\biggl(\prod_{j=1}^d \bigl(\lambda^{X_{1,j}}_{k_j}/\lambda^{X_{1,j}}_1\bigr)<C t \biggr)\\
&\leqslant&\sum\limits_{k=1}^\infty\hat\lambda^{X_\infty}_k\,\id\bigl(\hat\lambda^{X_\infty}_k< Ct \bigr).
\end{eqnarray*}
Since, by representation \eqref{conc_lambdaXd1_TPRE},  $C=\sum_{k=1}^{\infty}\hat\lambda^{X_\infty}_k<\infty$, we obtain the equality
\begin{eqnarray*}
\lim_{t\to0}\sup_{d\in\N}\sum\limits_{k=1}^\infty \nlambda^{X_d}_k\,\id\bigl(\nlambda^{X_d}_k<t \bigr)=\lim_{t\to0}\sum\limits_{k=1}^\infty\hat\lambda^{X_\infty}_k\,\id\bigl(\hat\lambda^{X_\infty}_k< Ct \bigr)=0,
\end{eqnarray*}
which is sufficient for $(i)$ by Proposition \ref{pr_nXdebounded_GRE}.\quad $\Box$

\subsection{Logarithmic asymptotics of the approximation complexity}\newparagraph
Here we will obtain criteria for the asymptotics \eqref{conc_nXde_LogAsymp_form} under the following assumption
\begin{eqnarray}\label{cond_UN_TPRE}
\lim_{d\to\infty}\max_{j=1,\ldots,d}\sum_{k=1}^\infty \nlambda^{X_{1,j}}_k\,\id\bigl(\nlambda^{X_{1,j}}_k< e^{-\tau b_d}\bigr)=0\quad\text{for all}\quad \tau>0,
\end{eqnarray}
which is rather weak, because of the next assertion.
\begin{Remark}\label{rem_Suff_UN_TPRE}
The condition $\lim\limits_{j\to\infty}\nlambda^{X_{1,j}}_1=1$ is sufficient for \eqref{cond_UN_TPRE} for any sequence $(b_d)_{d\in\N}$ such that $b_d\to\infty$, $d\to\infty$.
\end{Remark}
Indeed, let us fix $\tau>0$ and choose $j_\delta\in\N$ such that
\begin{eqnarray}\label{cond_nlambdaX1j1_Dominat1_TPRE}
1-\nlambda^{X_{1,j}}_1=\sum_{k=2}^\infty \nlambda^{X_{1,j}}_k\leqslant \delta\quad\text{for all}\quad j\geqslant j_\delta.
\end{eqnarray}
For all sufficiently large $d$ we get
\begin{eqnarray}\label{cond_nlambdaX1j1_Dominat2_TPRE}
\max_{j=1,\ldots,j_\delta} \sum_{k=1}^\infty \nlambda^{X_{1,j}}_k\id\bigl(\nlambda^{X_{1,j}}_k< e^{-\tau b_d}\bigr)\leqslant\delta.
\end{eqnarray}
Also for all $d$ such that $e^{-\tau b_d}<1-\delta$ we have
\begin{eqnarray*}
\max_{j=j_\delta,\ldots,d} \sum_{k=1}^\infty \nlambda^{X_{1,j}}_k\id\bigl(\nlambda^{X_{1,j}}_k< e^{-\tau b_d}\bigr)
&\leqslant& \max_{j=j_\delta,\ldots,d} \sum_{k=1}^\infty \nlambda^{X_{1,j}}_k\id\bigl(\nlambda^{X_{1,j}}_k/\nlambda^{X_{1,j}}_1< e^{-\tau b_d}/(1-\delta)\bigr)\\
&\leqslant& \max_{j=j_\delta,\ldots,d} \sum_{k=1}^\infty \nlambda^{X_{1,j}}_k\id\bigl(\nlambda^{X_{1,j}}_k/\nlambda^{X_{1,j}}_1< 1\bigr)\\
&\leqslant& \max_{j=j_\delta,\ldots,d}\sum_{k=2}^\infty \nlambda^{X_{1,j}}_k.
\end{eqnarray*}
Hence the condition \eqref{cond_UN_TPRE} follows from \eqref{cond_nlambdaX1j1_Dominat1_TPRE} and \eqref{cond_nlambdaX1j1_Dominat2_TPRE}.

Our next theorem gives the description of all possible functions $q$, which may appear  in $\eqref{conc_nXde_LogAsymp_form}$ under the assumption \eqref{cond_UN_TPRE}. Necessary notions and facts from probability theory can be found in Appendix.

\begin{Theorem}\label{th_nXde_LimitDistr_TPRE}
Let $(a_d)_{d\in\N}$ be a sequence, $(b_d)_{d\in\N}$ be a positive sequence such that $b_d\to\infty$, $d\to\infty$. Let a non-increasing function $q\colon (0,1)\to\R$ and a distribution function $G$ satisfy the equation   $q(\e)=G^{-1}(1-\e^2)$ for all $\e\in\ContSet(q)$. Suppose that under the condition \eqref{cond_UN_TPRE}, the following asymptotics holds
\begin{eqnarray}\label{th_nXde_LimitDistr_TPRE_cond}
\ln n^{X_d}(\e)=a_d+q(\e)b_d+o(b_d),\quad d\to\infty,\quad\text{for all}\quad \e\in\ContSet(q).
\end{eqnarray}
Then $G$ is self-decomposable with zero L\'evy spectral function on $(-\infty,0)$. If $G$ is non-degenerate, then $b_{d+1}/b_d\to1$  as $d\to\infty$.
\end{Theorem}

For proving this and next theorems we will introduce auxiliary random variables $U_j$, $j\in\N$, with the following distributions
\begin{eqnarray}\label{def_Uj}
\Probab\bigl(U_j=\bigl|\ln\nlambda^{X_{1,j}}_k\bigr|\bigr)=\mult^{X_{1,j}} \bigl(\nlambda^{X_{1,j}}_k\bigr)\cdot\nlambda^{X_{1,j}}_k,
\end{eqnarray}
where $j\in\N$ and $k\in\N$ such that $\lambda^{X_{1,j}}_k>0$.
\begin{Lemma}\label{lm_Uj}
Let $U_j$, $j\in\N$, be independent random variables with distributions $\eqref{def_Uj}$. Then we have
\begin{eqnarray}\label{lm_Uj_conc}
\sum_{k=1}^{\infty}\nlambda^{X_d}_k\,\id\bigl(\nlambda^{X_d}_k\geqslant  e^{-x}\bigr)=\Probab\Bigl(\sum_{j=1}^d U_j\leqslant x\Bigr)\quad\text{for all}\quad d\in\N,\quad x\in\R.
\end{eqnarray}
\end{Lemma}
\textbf{Proof of Lemma \ref{lm_Uj}.}\quad Let us fix $d\in\N$ and $x\in\R$. If $x<0$ then, by $\nlambda^{X_d}_k\leqslant 1$ and non-negativity of all $U_j$, $j\in\N$, both sides of \eqref{lm_Uj_conc} equal zero and hence it holds. Let $x\geqslant 0$. According to the multiplicative structure of $\nlambda^{X_d}_k$, $k\in\N$, we can write
\begin{eqnarray*}
\sum_{k=1}^{\infty}\nlambda^{X_d}_k\id\bigl(\nlambda^{X_d}_k\geqslant e^{-x}\bigr)
&=&\sum_{k\in\N^d} \prod_{j=1}^d\nlambda^{X_{1,j}}_{k_j}\id\biggr(\prod_{j=1}^d\nlambda^{X_{1,j}}_{k_j}\geqslant e^{-x}\biggr)\\
&=&\sum_{k\in\N^d} \prod_{j=1}^d\nlambda^{X_{1,j}}_{k_j}\id\Bigr(\sum_{j=1}^d\bigl|\ln\nlambda^{X_{1,j}}_{k_j}\bigr|\leqslant x\Bigr)\\
&=&\sum_{\substack{j=1,\ldots, d\\\mu_j\in\mathcal V^{X_{1,j}}}} \prod_{j=1}^d\mult^{X_{1,j}}(\mu_j)\,\mu_j\,\id\Bigr(\sum_{j=1}^d \bigl|\ln\mu_j\bigr|\leqslant x\Bigr),
\end{eqnarray*}
where $\mathcal V^{X_{1,j}}\colonequals \bigl\{\nlambda^{X_{1,j}}_k : k\in\N\bigr\}$ is the range of the sequence $(\nlambda^{X_{1,j}}_k)_{k\in\N}$ for every $j\in\N$. By \eqref{def_Uj}, we have
\begin{eqnarray*}
\sum_{k=1}^{\infty}\nlambda^{X_d}_k\,\id\bigl(\nlambda^{X_d}_k\geqslant e^{-x}\bigr)
&=&\sum_{\substack{ j=1,\ldots, d\\\mu_j\in\mathcal V^{X_{1,j}}}} \prod_{j=1}^d\Probab\bigl(U_j=\bigl|\ln \mu_j\bigr|\bigr)\id\Bigl(\sum_{j=1}^d\bigl|\ln \mu_j\bigr|\leqslant x\Bigr)\\
&=&\sum_{\substack{j=1,\ldots, d\\\mu_j\in\mathcal V^{X_{1,j}}}}\Probab\Bigl(U_j=\bigl|\ln\mu_j\bigr|,\, j\in\{1,\ldots,d\}\Bigr)\id\Bigl(\sum_{j=1}^d\bigl|\ln \mu_j\bigr|\leqslant x\Bigr)\\
&=&\Probab \Bigl(\sum_{j=1}^d U_j\leqslant x\Bigr).\quad\Box
\end{eqnarray*}

The described probabilistic construction was proposed by M. A. Lifshits and E. V. Tulyakova in the paper \cite{LifTul} in the context of approximation of tensor degrees of random elements, which have covariance operators with eigenvalues of unit multiplicity. We extend this approach to general tensor product-type random elements without any assumptions on marginal eigenvalue multiplicities.\\

\noindent\textbf{Proof of Theorem \ref{th_nXde_LimitDistr_TPRE}.}\quad According to Theorem \ref{th_nXde_LogAsymp_GRE}, the condition \eqref{th_nXde_LimitDistr_TPRE_cond} is equivalent to the convergence
\begin{eqnarray}\label{th_nXde_LogAsymp_TPRE_conv}
\lim_{d\to\infty} G^{X_d}_{a_d,b_d}(x)= G(x)\quad\text{for all}\quad x\in \ContSet(G),
\end{eqnarray}
where the functions $G^{X_d}_{a_d,b_d}$, $d\in\N$, are defined by \eqref{def_GXdadbd_GRE} for given  $a_d$, $b_d$, and $X_d$, $d\in\N$. Using Lemma~\ref{lm_Uj}, we can write 
\begin{eqnarray}\label{conc_GXdadbdSumUj_TPRE}
G^{X_d}_{a_d,b_d}(x)=\Probab \Biggl(\dfrac{\sum_{j=1}^d U_j - a_d}{b_d} \leqslant x\Biggr), \quad d\in\N, \quad x\in\R,
\end{eqnarray} 
where the independent random variables $U_j$, $j\in\N$, are distributed according to \eqref{def_Uj} with the tails
\begin{eqnarray*}
	\Probab\bigl(U_j>x\bigr)=\sum_{k=1}^\infty \nlambda^{X_{1,j}}_k\,\id\bigl(\nlambda^{X_{1,j}}_k< e^{-x}\bigr), \quad x\geqslant 0.
\end{eqnarray*}
From \eqref{cond_UN_TPRE} and non-negativity of $U_j$, $j\in\N$,  we have
\begin{eqnarray}\label{cond_UN_Uj_TPRE}
\lim_{d\to\infty}\max_{j=1,\ldots,d}\Probab\bigl(|U_j|>\tau b_d \bigr)=0\quad\text{for all}\quad \tau>0.
\end{eqnarray}
By Theorem \ref{th_LimitSelfDec_Probab} (see Appendix), the weak limit of $(G^{X_d}_{a_d,b_d})_{d\in\N}$, the function $G$, is necessarily self-decomposable. Let $L$ denote the  L\'evy spectral function of $G$. From $\condAtwo$ of Theorem \ref{th_NessSuffCond_ConvtoSelfDec_Probab} and from non-negativity of $U_j$, $j\in\N$, we get $L(x)=0$, $x<0$. The assertion about $(b_d)_{d\in\N}$  follows directly from Theorem \ref{th_LimitSelfDec_Probab}.\quad $\Box$\\

The next theorem provides a criterion for the asymptotics \eqref{th_nXde_LimitDistr_TPRE_cond}, where $q$ is a quantile of self-decomposable law.
\begin{Theorem}\label{th_nXde_LogAsymp_TPRE}
Let $(a_d)_{d\in\N}$ be a sequence, $(b_d)_{d\in\N}$ be a positive sequence such that $b_d\to\infty$, $d\to\infty$. Let a distribution function $G$ is self-decomposable with triplet $(\gamma, \sigma^2, L)$, where $L(x)=0$, $x<0$. Let a non-increasing function $q\colon (0,1)\to\R$ satisfy the equation   $q(\e)=G^{-1}(1-\e^2)$ for all $\e\in(0,1)$. Under the condition \eqref{cond_UN_TPRE}, for the asymptotics
\begin{eqnarray}\label{th_nXde_LogAsymp_TPRE_nXde}
\ln n^{X_d}(\e)=a_d+q(\e)b_d+o(b_d),\quad d\to\infty,\quad\text{for all}\quad \e\in(0,1),
\end{eqnarray}
the following ensemble of conditions is necessary and sufficient:
\begin{eqnarray*}
&&\text{$\condA$}\qquad\lim_{d\to\infty}\sum_{j=1}^d \sum_{k=1}^N \nlambda^{X_{1,j}}_k\id\bigl(\nlambda^{X_{1,j}}_k< e^{-x b_d}\bigr)=-L(x)\quad\text{for all}\quad x>0;\\
&&\text{$\condB$}\qquad\lim_{d\to\infty}\dfrac{1}{b_d}\Bigl(\sum_{j=1}^d M^{X_{1,j}}_{1,N}(\tau b_d)-a_d\Bigr)=\gamma+\int\limits_0^\tau \dfrac{x^3\,\dd L(x)}{1+x^2}-\int\limits_\tau^\infty \dfrac{x\,\dd L(x)}{1+x^2}\quad\text{for all}\quad \tau>0;\\
&&\text{$\condC$}\qquad\lim_{\tau\to0}\varlimsup_{d\to\infty}  \dfrac{1}{b_d^2}\sum_{j=1}^{d} \Bigl(M^{X_{1,j}}_{2,N}(\tau b_d)-\bigl(M^{X_{1,j}}_{1,N}(\tau b_d)\bigr)^2\Bigr)=\lim_{\tau\to0}\varliminf_{d\to\infty}\,\ldots\, =\sigma^2,
\end{eqnarray*}
where
\begin{eqnarray*}
M^{X_{1,j}}_{p,N}(x)\colonequals \sum_{k=1}^N \bigl|\ln\nlambda^{X_{1,j}}_k \bigr|^p\,\nlambda^{X_{1,j}}_k \, \id\bigl(\nlambda^{X_{1,j}}_k \geqslant e^{-x}\bigr),\quad x\geqslant0,\quad p\in\{1,2\},
\end{eqnarray*}
and $N$ is $\infty$ or any natural number such that
\begin{eqnarray}\label{th_nXde_LogAsymp_TPRE_cond_N}
\sum_{j=1}^{\infty} \sum_{k=N+1}^\infty\nlambda^{X_{1,j}}_k<\infty.
\end{eqnarray}
\end{Theorem}
\textbf{Proof of Theorem \ref{th_nXde_LogAsymp_TPRE}.}\quad We first show that $\ContSet(q)=(0,1)$. If $G$ is degenerate, then $q$ is constant on $(0,1)$ and, consequently, continuous. For non-degenerate case the assertion follows from Remarks \ref{rem_Contq_GRE} and  \ref{rem_StrIncSelfDec_Probab} (see Appendix). Thus, by Theorem \ref{th_nXde_LogAsymp_GRE}, the condition \eqref{th_nXde_LogAsymp_TPRE_nXde} is equivalent to the convergence \eqref{th_nXde_LogAsymp_TPRE_conv}, where $G^{X_d}_{a_d,b_d}$, $d\in\N$, are defined by \eqref{def_GXdadbd_GRE} for given  $a_d$, $b_d$, and $X_d$, $d\in\N$. According to \eqref{conc_GXdadbdSumUj_TPRE}, $G^{X_d}_{a_d,b_d}$, $d\in\N$, are distribution functions of centered and normalized  sums of non-negative independent random variables $U_j$, $j\in\N$, satisfying \eqref{cond_UN_Uj_TPRE}. Consequenly, for the convergence \eqref{th_nXde_LogAsymp_TPRE_conv} the conditions $\condAone$, $\condB$, and $\condC$ of Theorem \ref{th_NessSuffCond_ConvtoSelfDec_Probab} (see Appendix) are necessary and sufficient (we set $Y_j\colonequals U_j$, $j\in\N$, and $A_d\colonequals a_d$, $B_d\colonequals b_d$, $d\in\N$ in Theorem \ref{th_NessSuffCond_ConvtoSelfDec_Probab}). For the case $N=\infty$ it directly yields the conditions $\condA$, $\condB$, and $\condC$ of the theorem. Indeed, it is easily seen that
\begin{eqnarray*}
\Expec \Bigl[U_j^p\,\id\bigl(|U_j|\leqslant x\bigr)\Bigr]&=&
M^{X_{1,j}}_{p,\infty}(x),\quad x\geqslant0,\quad  p\in\{1,2\}.
\end{eqnarray*}
Hence for any $x\geqslant 0$
\begin{eqnarray*}
\Disp  \Bigl[ U_j\,\id\bigl(|U_j|\leqslant x\bigr)\Bigr]
&=&\Expec \bigl[\,\ldots\,\bigr]^2- \bigl(\Expec  \bigl[\,\ldots\,\bigr]\bigr)^2\\
&=&M^{X_{1,j}}_{2,\infty}(x)-\bigl(M^{X_{1,j}}_{1,\infty}(x)\bigr)^2.
\end{eqnarray*}

Suppose that there exists a natural number $N=N'$ that satisfies \eqref{th_nXde_LogAsymp_TPRE_cond_N}. Now we  show that the conditions $\condA$, $\condB$, and $\condC$  for $N=\infty$ are respectively equivalent to the same ones for $N=N'$.

For arbitrary small $\delta>0$ we can find $j_\delta\in\N$ that yields $\sum_{j>j_\delta} \sum_{k>N'}\nlambda^{X_{1,j}}_k\leqslant \delta$. Thus for any $x>0$ we have
\begin{eqnarray*}
\sum_{j=1}^d \sum_{k>N'}\nlambda^{X_{1,j}}_k \id\bigl(\nlambda^{X_{1,j}}_k<e^{-x b_d} \bigr)
&\leqslant& \sum_{j\leqslant j_\delta}\sum_{k>N'}\nlambda^{X_{1,j}}_k \id\bigl(\nlambda^{X_{1,j}}_k<e^{-x b_d} \bigr)
+\sum_{j>j_\delta} \sum_{k>N'}\nlambda^{X_{1,j}}_k\\
&\leqslant& j_\delta\cdot \max_{j=1,\ldots,j_\delta}\sum_{k=1}^\infty\nlambda^{X_{1,j}}_k \id\bigl(\nlambda^{X_{1,j}}_k<e^{-x b_d} \bigr)+\delta.
\end{eqnarray*}
By \eqref{cond_UN_TPRE}, the last expression can be made arbitrary small by the choice of sufficiently small $\delta$ and large $d$. This proves the equivalence of 
$\condA$ for $N=\infty$ and $\condA$ for $N=N'$ under the condition \eqref{th_nXde_LogAsymp_TPRE_cond_N} with $N=N'$.

To continue the proof, we need the following relation
\begin{eqnarray}\label{conc_UN1_TPRE}
\lim_{d\to\infty}\dfrac{1}{b_d}\,\max_{j=1,\ldots,d} \sum_{k=1}^\infty \bigl|\ln\nlambda^{X_{1,j}}_k \bigr|\,\nlambda^{X_{1,j}}_k \, \id\bigl(\nlambda^{X_{1,j}}_k \geqslant e^{-\tau b_d}\bigr)=0\quad\text{for all}\quad\tau>0, 
\end{eqnarray}
which follows from \eqref{cond_UN_TPRE} and the inequalities
\begin{eqnarray*}
\dfrac{1}{b_d}\,\max_{j=1,\ldots,d} \sum_{k=1}^\infty \bigl|\ln\nlambda^{X_{1,j}}_k \bigr|\,\nlambda^{X_{1,j}}_k \, \id\bigl(\nlambda^{X_{1,j}}_k \geqslant e^{-\tau b_d}\bigr)
&\leqslant&
\dfrac{1}{b_d}\,\max_{j=1,\ldots,d} \biggl\{\sum_{k=1}^\infty \bigl|\ln\nlambda^{X_{1,j}}_k \bigr|\,\nlambda^{X_{1,j}}_k \, \id\bigl(\nlambda^{X_{1,j}}_k \geqslant e^{-h b_d}\bigr)\\
&&{}+\sum_{k=1}^\infty \bigl|\ln\nlambda^{X_{1,j}}_k \bigr|\,\nlambda^{X_{1,j}}_k \, \id\Bigl(\nlambda^{X_{1,j}}_k \in\bigl[e^{-\tau b_d}, e^{-h b_d}\bigl)\Bigr)\biggl\}\\
&\leqslant& h+\tau \max_{j=1,\ldots,d}\sum_{k=1}^\infty\nlambda^{X_{1,j}}_k \id\bigl(\nlambda^{X_{1,j}}_k<e^{-h b_d} \bigr),
\end{eqnarray*}
where $h<\tau$ is arbitrary small positive number. 

Next, for any $\tau>0$ we estimate the following sums
\begin{eqnarray*}
\dfrac{1}{b_d}\sum_{j=1}^d \sum_{k>N'} \bigl|\ln\nlambda^{X_{1,j}}_k \bigr|\,\nlambda^{X_{1,j}}_k \, \id\bigl(\nlambda^{X_{1,j}}_k \geqslant e^{-\tau b_d}\bigr)
&\leqslant& \dfrac{1}{b_d}\sum_{j\leqslant j_\delta} \sum_{k>N'} \bigl|\ln\nlambda^{X_{1,j}}_k \bigr|\,\nlambda^{X_{1,j}}_k \, \id\bigl(\nlambda^{X_{1,j}}_k \geqslant e^{-\tau b_d}\bigr)\\
&&{}+\dfrac{1}{b_d}\sum_{j>j_\delta} \sum_{k>N'}\bigl|\ln\nlambda^{X_{1,j}}_k \bigr|\,\nlambda^{X_{1,j}}_k \, \id\bigl(\nlambda^{X_{1,j}}_k \geqslant e^{-\tau b_d}\bigr)\\
&\leqslant& \dfrac{j_\delta}{b_d}\,\max_{j=1,\ldots,j_\delta} \sum_{k=1}^\infty \bigl|\ln\nlambda^{X_{1,j}}_k \bigr|\,\nlambda^{X_{1,j}}_k \, \id\bigl(\nlambda^{X_{1,j}}_k \geqslant e^{-\tau b_d}\bigr)\\
&&+\tau \sum_{j>j_\delta}\sum_{k>N'}\nlambda^{X_{1,j}}_k\\
&\leqslant& \dfrac{j_\delta}{b_d}\,\max_{j=1,\ldots,d} \sum_{k=1}^\infty \bigl|\ln\nlambda^{X_{1,j}}_k \bigr|\,\nlambda^{X_{1,j}}_k \, \id\bigl(\nlambda^{X_{1,j}}_k \geqslant e^{-\tau b_d}\bigr)+\tau \delta.
\end{eqnarray*}
Here the last expression can be made arbitrary small by the choice of sufficiently small $\delta$ and large $d$ in view of \eqref{conc_UN1_TPRE}. Thus we have the equivalence of 
$\condB$ for $N=\infty$ and $\condB$ for $N=N'$ under the condition \eqref{th_nXde_LogAsymp_TPRE_cond_N} with $N=N'$.

To obtain the equivalence $\condC$ for $N=\infty$ and $\condC$ for $N=N'$, it is sufficient to show that
\begin{eqnarray*}
\lim_{d\to\infty}\dfrac{1}{b_d^2}\sum_{j=1}^{d}\sum_{k>N'} \bigl|\ln\nlambda^{X_{1,j}}_k \bigr|^2\,\nlambda^{X_{1,j}}_k \, \id\bigl(\nlambda^{X_{1,j}}_k \geqslant e^{-\tau b_d}\bigr)=0\quad\text{for all}\quad \tau>0.
\end{eqnarray*}
This follows from the estimate
\begin{eqnarray*}
\dfrac{1}{b_d^2}\sum_{k>N'} \bigl|\ln\nlambda^{X_{1,j}}_k \bigr|^2\,\nlambda^{X_{1,j}}_k \, \id\bigl(\nlambda^{X_{1,j}}_k \geqslant e^{-\tau b_d}\bigr)\leqslant \dfrac{\tau}{b_d}\sum_{k>N'} \bigl|\ln\nlambda^{X_{1,j}}_k \bigr|\,\nlambda^{X_{1,j}}_k \, \id\bigl(\nlambda^{X_{1,j}}_k \geqslant e^{-\tau b_d}\bigr)
\end{eqnarray*}
and the previous conclusions. \quad $\Box$

\subsection{Applications to tensor products of Euler integrated processes}\newparagraph
Let us consider Gaussian random process $E_r(t)$, $t\in[0,1]$, with zero mean and the following correlation function
\begin{eqnarray*}
\CorrFunc^{E_r}(t,s)\colonequals\int\limits_{[0,1]^{r}} \min\{t,x_1\} \min\{x_1,x_2\}\ldots  \min\{x_{r},s\} \dd x_1\ldots \dd x_{r}, 
\end{eqnarray*}
where $t, s\in[0,1]$ and $r$ is a non-negative integer. Such process is called the \textit{Euler integrated process}. It is connected with the standard Wiener process $W(t)$, $t\in[0,1]$, by the following integration scheme:
\begin{eqnarray*}
E_r(t)=(-1)^{a_1+\ldots+a_r} \int\limits_{a_r}^{t} \int\limits_{a_{r-1}}^{t_{r-1}}\ldots \int\limits_{0}^{t_2}\int\limits_{1}^{t_1}  W(s) \dd s\, \dd t_1 \ldots \dd t_{r-1}, \quad t\in[0,1],
\end{eqnarray*}
where $a_{2k-1}:=1$, $a_{2k}:=0$, $k\in\N$. The Euler integrated process is well known object: related boundary value problems were considered in \cite{ChangHa}, small ball probabilities were investigated in the papers \cite{GaoHanTor} and \cite{NazNik}.

Let us consider the sequence of the Euler integrated processes $E_{r_j}(t)$, $t\in[0,1]$, $j\in\N$, with correlation functions $\CorrFunc^{E_{r_j}}$, $j\in\N$, respectively. 
By $(r_j)_{j\in\N}$ we will always mean a non-decreasing sequence of non-negative integers. Consider a sequence of zero-mean random fields $\EulerField_d(t)$, $t\in[0,1]^d$, $d\in\N$, with the following correlation functions
\begin{eqnarray*}
\CorrFunc^{\EulerField_d}(t,s)=\prod\limits_{j=1}^d \CorrFunc^{E_{r_j}} (t_j,s_j),\quad t,s\in[0,1]^d,\quad d\in\N.
\end{eqnarray*}
We consider every process $E_{r_j}(t)$, $t\in[0,1]$, as a random element $E_{r_j}$ of the space $L_2([0,1])$. The covariance operator $K^{E_{r_j}}$ of $E_{r_j}$ is the integration operator with kernel $\CorrFunc^{E_{r_j}}$. The eigenpairs of $K^{E_{r_j}}$ are exactly known (see \cite{ChangHa} and \cite{GaoHanTor}):
\begin{eqnarray*}
\lambda^{E_{r_j}}_k=\dfrac{1}{\bigl(\pi(k-1/2)\bigr)^{2r_j+2}},\quad \psi^{E_{r_j}}_k(t)=\sqrt{2} \sin\bigl(\pi(k-1/2)t\bigr),\quad k\in\N, \quad t\in[0,1].
\end{eqnarray*}
It is easily seen that every field $\EulerField_d(t)$, $t\in[0,1]^d$, as a random element $\EulerField_d$ of $L_2([0,1]^d)$, has a covariance operator of the form $K^{\EulerField_d}=\otimes_{j=1}^d K^{E_{r_j}}$, with eigenpairs \eqref{conc_eigenpairs_TPRE} (we set $X_d=\EulerField_d$, $X_{1,j}=E_{r_j}$). Therefore, by definition from Section 2, $\EulerField_d$ is a tensor product of $E_{r_j}$, i.e. $\EulerField_d=\otimes_{j=1}^d E_{r_j}$. 

Let us consider the sequence $\APP$ of approximation problems for $\EulerField_d$, $d\in\N$.  The criteria of all types of tractability for $\APP$ were obtained in the paper \cite{LifPapWoz2}. We recall a part of this result here.
\begin{Theorem}\label{th_Tract_Euler}
$\APP$ is \textit{weakly tractable} iff 
\begin{eqnarray} \label{th_Tract_Euler_cond_WT}
\lim_{j\to\infty} r_j=\infty;
\end{eqnarray}
$\APP$ is \textit{quasi-polynomially tractable} iff
\begin{eqnarray}\label{th_Tract_Euler_cond_QPT}
\sup_{d\in\N} \dfrac{1}{\ln_+ d}\sum_{j=1}^d (r_j+1)\, 3^{-2r_j-2}<\infty;
\end{eqnarray}
$\APP$ is \textit{polynomially tractable} iff it is \textit{strongly polynomially tractable} iff
\begin{eqnarray}\label{th_Tract_Euler_cond_SPT}
\sum_{j=1}^{\infty} 3^{-2\tau (r_j+1)}<\infty\quad\text{for some}\quad \tau\in(0,1).
\end{eqnarray}
\end{Theorem}
In fact, depending on the type of tractability, the previous theorem provides a majorant (see \eqref{def_WT}--\eqref{def_SPT}) for the approximation complexity $n^{\EulerField_d}(\e)$ for all $\e\in(0,1)$ and $d\in\N$.

We will investigate $n^{\EulerField_d}(\e)$ for arbitrarily fixed $\e\in(0,1)$ and $d\to\infty$. In order to compare our results with tractability bounds, we restrict ourselves only to quasi-polinomially tractable sequences $\APP$ with the following behaviour of $(r_j)_{j\in\N}$:
\begin{eqnarray}\label{cond_lambdaErj2_Asymp_Euler}
3^{-2r_j-2}\sim \dfrac{\beta}{ j (\ln j)^p},\quad j\to\infty,
\end{eqnarray}
where $\beta>0$ and $p\geqslant 1$.

\begin{Remark}
Under the assumption \eqref{cond_lambdaErj2_Asymp_Euler} for some $\beta>0$ and $p\in\R$, $\APP$ is weakly tractable, but not strongly polynomially tractable. $\APP$ is quasi-polynomially tractable iff $p\geqslant1$. 
\end{Remark}
Indeed, weak tractability immediately follows from \eqref{th_Tract_Euler_cond_WT}. It is easy to check that \eqref{th_Tract_Euler_cond_SPT} does not hold. Consider the behaviour of sums from \eqref{th_Tract_Euler_cond_QPT} as $d\to\infty$
\begin{eqnarray*}
\dfrac{1}{\ln_+ d}\sum_{j=1}^d (1+r_j)\, 3^{-2r_j-2}
&\sim&\dfrac{\beta}{(2\ln 3)\ln  d}\,\sum_{j=1}^d \dfrac{(\ln j)^{1-p}}{j}\\
&\sim&\dfrac{\beta}{(2\ln 3)\ln  d}\int\limits_2^d\dfrac{(\ln t)^{1-p}}{t} \, \dd t\\
&\sim& \dfrac{\beta}{(2\ln 3)\ln  d}\cdot\biggl(\dfrac{(\ln d)^{2-p}}{2-p}\,\id(p\ne 2)+ (\ln\ln d)\,\id(p=2)\biggr) \\
&=&\dfrac{\beta}{2\ln 3}\cdot\biggl(\dfrac{(\ln d)^{1-p}}{2-p}\,\id(p\ne 2)+ \dfrac{\ln\ln d}{\ln d}\,\id(p=2)\biggr).
\end{eqnarray*} 
From this we conclude the second assertion of the remark.

\begin{Proposition}\label{pr_nXde_Boundedness_Euler}
Suppose that \eqref{cond_lambdaErj2_Asymp_Euler} holds for some $\beta>0$ and $p\geqslant 1$. If $p>1$ then
\begin{eqnarray*}
\sup_{d\in\N} n^{\EulerField_d}(\e)<\infty\quad\text{for all}\quad \e\in(0,1).
\end{eqnarray*}
If $p=1$ then 
\begin{eqnarray*}
\lim\limits_{d\to\infty}n^{\EulerField_d}(\e)=\infty\quad\text{for all}\quad\e\in(0,1).
\end{eqnarray*}
\end{Proposition}
\textbf{Proof of Proposition \ref{pr_nXde_Boundedness_Euler}.}\quad Consider the following sum 
\begin{eqnarray*}
\sum_{j=1}^{d} \sum_{k=2}^{\infty} \dfrac{\lambda^{E_{r_j}}_k}{\lambda^{E_{r_j}}_1} =\sum_{j=1}^{d} \sum_{k=2}^{\infty}\dfrac{1}{(2k-1)^{2r_j+2}}=\sum_{j=1}^{d} c(r_j) 3^{-2r_j-2},
\end{eqnarray*}
where $c(r_j)\colonequals 1+\sum_{k=3}^{\infty}\bigl((2k-1)/3\bigr)^{-2r_j-2}$. Under the assumption \eqref{cond_lambdaErj2_Asymp_Euler} for some $\beta>0$ and $p\geqslant1$, $c(r_j)\to1$, $j\to\infty$, and the series $\sum_{j=1}^{\infty} 3^{-2r_j-2}$ converges for $p>1$ and it diverges for $p=1$. Applying Proposition \ref{pr_nXdebounded_TPRE}, we obtain the required assertions.\quad $\Box$\\

From this proposition we can see that for the cases with $p>1$ the bound \eqref{def_QPT} for $n^{X_d}(\e)$ can be rather crude under the setting ``$\e$ is fixed, $d\to\infty$'', because it depends on $d$, whereas $n^{X_d}(\e)$ is bounded on $d$.

For the case $p=1$ we can find logarithmic asymptotics of the approximation complexity. Here the \textit{convolution powers of the Dickman distribution} appear (see Appendix).
\begin{Proposition}\label{pr_nXde_Dickman_Euler}
Under the assumption \eqref{cond_lambdaErj2_Asymp_Euler} with $p=1$ and $\beta>0$, we have
\begin{eqnarray}\label{pr_nXde_Dickman_Euler_conc}
\ln n^{\EulerField_d}(\e)=D_\beta^{-1}(1-\e^2)\ln d+o(\ln d),\quad d\to\infty,\quad\text{for all}\quad\e\in(0,1),
\end{eqnarray}
where $D_\beta$ is the distribution function of $\beta$-convolution power of the Dickman distribution.
\end{Proposition}
\textbf{Proof of Proposition \ref{pr_nXde_Dickman_Euler}.}\quad Let us write the expression for traces of $K^{E_{r_j}}$ in the following form:
\begin{eqnarray*}
\Lambda^{E_{r_j}}=\sum_{k=1}^{\infty}\lambda^{E_{r_j}}_k=\sum_{k=1}^{\infty} \bigl(\pi(k-1/2)\bigr)^{-2r_j-2}=\dfrac{1}{\omega_{r_j}}\cdot \biggl(\dfrac{2}{\pi}\biggr)^{2r_j+2},
\end{eqnarray*}
where $\omega_{r_j}\colonequals \bigl(\sum_{k=1}^{\infty} (2k-1)^{-2r_j-2}\bigr)^{-1}$, $j\in\N$. According to the notation, we have
\begin{eqnarray*}
\nlambda^{E_{r_j}}_k=\dfrac{\omega_{r_j}}{(2k-1)^{2r_j+2}},\quad k\in\N,\quad j\in\N.
\end{eqnarray*}

We see that $\nlambda^{E_{r_j}}_1=\omega_{r_j}\to 1$ as $j\to\infty$. Hence, by Remark \ref{rem_Suff_UN_TPRE}, the condition \eqref{cond_UN_TPRE} holds for any sequence $(b_d)_{d\in\N}$ such that $b_d\to\infty$, $d\to\infty$. Next, consider the sums
\begin{eqnarray*}
\sum_{j=1}^{d}\sum_{k=3}^\infty\nlambda^{E_{r_j}}_k=\sum_{j=1}^{d} \sum_{k=3}^{\infty} \dfrac{\omega_{r_j}}{(2k-1)^{2r_j+2}}\leqslant C\cdot \sum_{j=1}^{d}  5^{-2r_j-2},
\end{eqnarray*}
where we set $C_1:=1+\sum_{k=4}^{\infty}\bigl(5/(2k-1)\bigr)^{2}<\infty$. Under our assumptions on $(r_j)_{j\in\N}$, we get
\begin{eqnarray*}
5^{-2r_j-2}\sim \biggl(\dfrac{\beta}{j\, \ln j} \biggr)^{\tfrac{\ln 5}{\ln 3}},\quad j\to\infty,
\end{eqnarray*}
which gives the convergence of the series $\sum_{j=1}^{\infty}\sum_{k=3}^\infty\nlambda^{E_{r_j}}_k$. Thus, in order to obtain the asymptotics \eqref{pr_nXde_Dickman_Euler_conc}, it is sufficient to check conditions $\condA$--$\condC$ of Theorem \ref{th_nXde_LogAsymp_TPRE} with $X_d=\EulerField_d$, $X_{1,j}=E_{r_j}$, $G=D_\beta$, $a_d=0$, $b_d=\ln_+ d$, and $N=2$. Here $D_\beta$ has the following triplet $\gamma\colonequals\beta\pi/4$, $\sigma^2\colonequals 0$, $L(x)\colonequals\beta \ln x\, \id\bigl(x\in(0,1]\bigr)$ (see Appendix).

Let us check the condition $\condA$ of Theorem \ref{th_nXde_LogAsymp_TPRE}:
\begin{eqnarray}\label{cond_sumnlambdaErj12_Euler}
\lim_{d\to\infty}\sum_{j=1}^d \sum_{k=1}^2 \nlambda^{E_{r_j}}_k\id\bigl(\nlambda^{E_{r_j}}_k< e^{-x \ln d}\bigr)=-\beta \ln x\, \id\bigl(x\in(0,1]\bigr)\quad\text{for all}\quad x>0.
\end{eqnarray}
Since $\nlambda^{E_{r_j}}_1\to 1$ as $j\to\infty$, for any $x>0$ and sufficiently large $d$ we have $\min_{j\in\N}\nlambda^{E_{r_j}}_1>e^{-x\ln d}$ and, consequently,
\begin{eqnarray*}
\sum_{j=1}^d \sum_{k=1}^2 \nlambda^{E_{r_j}}_k\id\bigl(\nlambda^{E_{r_j}}_k< e^{-x \ln d}\bigr)=\sum_{j=1}^d \nlambda^{E_{r_j}}_2\id\bigl(\nlambda^{E_{r_j}}_2< e^{-x \ln d}\bigr).
\end{eqnarray*}
The right-hand side of this equation tends to $-\beta \ln x\, \id\bigl(x\in(0,1]\bigr)$ as $d\to\infty$ for any $x>0$, iff the sums $\sum_{j=1}^d  3^{-2r_j-2}\,\id\bigl(3^{2r_j+2}> e^{x \ln d}\bigr)$ tend to the same limit as $d\to\infty$ for any $x>0$. This follows from the equality $\nlambda^{E_{r_j}}_2=\omega_{r_j} 3^{-2r_j-2}$ with $\omega_{r_j}\to1$, $j\to\infty$, and from the continuity of the limit function. Thus \eqref{cond_sumnlambdaErj12_Euler} is equivalent to
\begin{eqnarray}\label{conc_sum3rj_ConvtoLevyFunc_Euler}
\lim_{d\to\infty}\sum_{j=1}^d  3^{-2r_j-2}\,\id\bigl(3^{2r_j+2}> e^{x \ln d}\bigr)= - \beta\ln x\, \id\bigl(x\in(0,1]\bigr)\quad \text{for all}\quad x>0.
\end{eqnarray}
For $x\geqslant1$ we have  $3^{2r_j+2}>e^{x \ln d}$, $j=1, \ldots,d$, for all sufficiently large $d$. Therefore \eqref{conc_sum3rj_ConvtoLevyFunc_Euler} obviously holds in this case.
For $x\in(0,1)$ we set $j_{d,x}:=\min\bigl\{j\in\N: 3^{2r_j+2}>e^{x\ln d}\bigr\}$. Using Lemmas \ref{lm_AsymInv_deBruijn_SVF} and \ref{lm_Bekessy_SVF} from Appendix, we have
\begin{eqnarray}\label{conc_jdx_asymp}
j_{d,x}\sim \beta d^x (\ln )^\#(\beta d^x)\sim\dfrac{\beta\,d^x}{x\ln d}, \quad d\to\infty.
\end{eqnarray}
Under our assumptions on $(r_j)_{j\in\N}$, we have
\begin{eqnarray*}
\sum_{j=1}^d  3^{-2r_j-2}\,\id\bigl(3^{2r_j+2}> e^{x \ln d}\bigr)= \sum_{j=j_{d,x}}^{d} 3^{-2r_j-2}\sim \sum_{j=j_{d,x}}^{d} \dfrac{\beta}{j \ln j}, \quad d\to\infty.
\end{eqnarray*}
Using \eqref{conc_jdx_asymp} and the asymptotics (see \cite{MakGolLodPod},  2.13, p. 21)
\begin{eqnarray}\label{conc_lnlnj_Asymp_Euler}
\sum_{k=1}^{n}\dfrac{1}{k \ln k}=\ln\ln n+c+o(1), \quad n\to\infty,
\end{eqnarray}
with some constant $c$, we obtain
\begin{eqnarray*}
\sum_{j=j_{d,x}}^{d} \dfrac{\beta}{ j \ln j}&=&\beta\bigl(\ln\ln d - \ln\ln j_{d,x}\bigr)+o(1)\\
&=&\beta\bigl(\ln\ln d - \ln\ln d^x \bigr)+o(1)\\
&=&-\beta\ln x+o(1),
\end{eqnarray*}
as $d\to\infty$. Thus we have the convergence \eqref{conc_sum3rj_ConvtoLevyFunc_Euler}. 

Next, we check the condition $\condB$ of Theorem \ref{th_nXde_LogAsymp_TPRE}: 
\begin{eqnarray*}
\lim_{d\to\infty}\dfrac{1}{\ln d}\sum_{j=1}^d \sum_{k=1}^2 \bigl|\ln\nlambda^{E_{r_j}}_k \bigr|\,\nlambda^{E_{r_j}}_k \, \id\bigl(\nlambda^{E_{r_j}}_k \geqslant e^{-\tau \ln d}\bigr)= \dfrac{\beta\pi}{4}+\int\limits_0^{\min\{\tau,1\}} \dfrac{\beta x^2\,\dd x}{1+x^2}-\int\limits_{\min\{\tau,1\}}^1 \dfrac{\beta\dd x}{1+x^2}
\end{eqnarray*}
for all $\tau>0$. Here the right-hand side exactly equals to $\beta \min\{\tau,1\}$. In the left-hand side we estimate
\begin{eqnarray*}
\sum_{j=1}^d \bigl|\ln\nlambda^{E_{r_j}}_1 \bigr|\,\nlambda^{E_{r_j}}_1
\leqslant \sum_{j=1}^d \bigl(1-\nlambda^{E_{r_j}}_1\bigr)=\sum_{j=1}^d \sum_{k=2}^{\infty} \dfrac{\omega_{r_j}}{(2k-1)^{2r_j+2}}\leqslant C_2 \sum_{j=1}^d 3^{-2r_j-2},
\end{eqnarray*}
where we set $C_2:=1+\sum_{k=3}^{\infty}\bigl(3/(2k-1)\bigr)^{2}<\infty$. Next, by \eqref{conc_lnlnj_Asymp_Euler}, we see that
\begin{eqnarray*}
\dfrac{C_2}{\ln d} \sum_{j=1}^d 3^{-2r_j-2}\sim \dfrac{C_2}{\ln d} \sum_{j=1}^d \dfrac{\beta}{j\ln j}\sim \dfrac{C_2\beta\ln\ln d}{\ln d}\to 0,\quad d\to\infty.
\end{eqnarray*}
Hence we only need to prove the convergence
\begin{eqnarray}\label{conc_sumlnnlambdaErj2_Euler}
\lim_{d\to\infty}\dfrac{1}{\ln d}\sum_{j=1}^d\bigl|\ln\nlambda^{E_{r_j}}_2 \bigr|\,\nlambda^{E_{r_j}}_2 \, \id\bigl(\nlambda^{E_{r_j}}_2 \geqslant e^{-\tau \ln d}\bigr)= \beta\min\{\tau,1\}\quad\text{for all}\quad \tau>0.
\end{eqnarray}

Since $\nlambda^{E_{r_j}}_2=\omega_{r_j} 3^{-2r_j-2}$ with $\omega_{r_j}\to1$, $j\to\infty$, and  the limit function $\tau\mapsto\beta\min\{\tau,1\}$ is continuous at any $\tau>0$, \eqref{conc_sumlnnlambdaErj2_Euler} is equivalent to
\begin{eqnarray*}
\lim_{d\to\infty}\dfrac{1}{\ln d}\sum_{\substack{j=1,\ldots,d\\j<j_{d,\tau}}}\ln \bigl(3^{2r_j+2}\bigr) \, 3^{-2r_j-2} = \beta\min\{\tau,1\}\quad\text{for all}\quad \tau>0.
\end{eqnarray*}
But this follows from \eqref{conc_jdx_asymp} and the following equivalences
\begin{eqnarray*}
\sum_{\substack{j=1,\ldots,d\\j<j_{d,\tau}}}\ln \bigl(3^{2r_j+2}\bigr) \, 3^{-2r_j-2}\sim \sum_{\substack{j=1,\ldots,d\\j<j_{d,\tau}}} \dfrac{\beta}{j}\sim \beta\ln\min\{d, j_{d,\tau}\}\sim \beta\min\{\tau,1\}\ln d,\quad d\to\infty.
\end{eqnarray*}

Now we check the condition $\condC$ of Theorem \ref{th_nXde_LogAsymp_TPRE}.  It is sufficient to prove
\begin{eqnarray*}
\lim_{\tau\to0}\varlimsup_{d\to\infty}\sum_{j=1}^d \sum_{k=1}^2 \dfrac{\bigl|\ln\nlambda^{E_{r_j}}_k \bigr|^2}{(\ln d)^2}\,\nlambda^{E_{r_j}}_k \, \id\bigl(\nlambda^{E_{r_j}}_k \geqslant e^{-\tau \ln d}\bigr)=0.
\end{eqnarray*}
But this follows from \eqref{conc_sumlnnlambdaErj2_Euler} and the next bound
\begin{eqnarray*}
\varlimsup_{d\to\infty}\sum_{j=1}^d \sum_{k=1}^2 \dfrac{\bigl|\ln\nlambda^{E_{r_j}}_k \bigr|^2}{(\ln d)^2}\,\nlambda^{E_{r_j}}_k \, \id\bigl(\nlambda^{E_{r_j}}_k \geqslant e^{-\tau \ln d}\bigr)
&\leqslant& \varlimsup_{d\to\infty}\sum_{j=1}^d \sum_{k=1}^2 \dfrac{\tau\bigl|\ln\nlambda^{E_{r_j}}_k \bigr|}{\ln d}\,\nlambda^{E_{r_j}}_k \, \id\bigl(\nlambda^{E_{r_j}}_k \geqslant e^{-\tau \ln d}\bigr)\\
&=&\beta\tau \min\{\tau,1\}.
\end{eqnarray*}

Thus we have checked the conditions $\condA$--$\condC$ of Theorem  \ref{th_nXde_LogAsymp_TPRE}. Hence, by this theorem we have the required asymptotics. \quad $\Box$

\section{Approximation of tensor degree-type random elements}\newparagraph
In this section we consider sequences of tensor degree-type random elements. Here we deal with an important particular case of the linear tensor approximation problems, which were considered in Section 4. Let $X$ be a zero-mean random element of a separable Hilbert space $H$ and assume that $\Expec \|X\|^2_H<\infty$. Let $X_d=X^{\otimes d}$ for every $d\in\N$. From \eqref{conc_eigenpairs_TPRE} we see that eigenvalues and eigenvectors of $K^{X_d}$ are respectively the following:
\begin{eqnarray*} 
\prod_{j=1}^d\lambda^{X}_{k_j},\quad \bigotimes_{j=1}^d\psi^{X}_{k_j},\qquad k_1, k_2, \ldots,k_d\in\N.
\end{eqnarray*}
By formula \eqref{conc_LambdaXd_TPRE}, for traces of $K^{X_d}$ we have $\Lambda^{X_d}=(\Lambda^{X})^d$, $d\in\N$. As in Section 4, we always assume that $\lambda^{X}_1>0$.

Let us proceed with some elementary remarks following from the definitions. If $\lambda^X_1=\Lambda^X$ then by formula \eqref{conc_nXde1} we have $n^{X_d}(\e)=1$ for all $\e\in(0,1)$ and $d\in\N$. In more interesting case, when $\lambda^X_1<\Lambda^X$, from the equality $\lambda^{X_d}_1=(\lambda^X_1)^d$ and lower bound \eqref{conc_nXde_lambdaXd1_ineq_GRE} we obtain the inequality
\begin{eqnarray}\label{conc_nXde_lambdaXd1_ineq_TDRE}
n^{X_d}(\e)\geqslant (1-\e^2)\,\bigl(\Lambda^X/\lambda^X_1\bigr)^d\quad\text{for all}\quad \e\in(0,1),\,\, d\in\N.
\end{eqnarray}
Therefore $n^{X_d}(\e)\to\infty$ as $d\to\infty$ at least exponentially for all fixed $\e\in(0,1)$ i.e. here we always have the curse of dimensionality.

In general, bound \eqref{conc_nXde_lambdaXd1_ineq_TDRE} can not be improved. Indeed, consider the case when
\begin{eqnarray}\label{cond_lambdaXk_triv_TDRE}
\exists\, l^X\in\N\qquad \lambda^X_k=\lambda^X_1\id(k\leqslant l^X),\quad k\in\N.
\end{eqnarray}
Under this condition, we have $l^X=\Lambda^X/\lambda^X_1$ and $\lambda^{X_d}_k=\lambda^{X_d}_1\id\bigl(k\leqslant (l^X)^d\bigr)$, $k\in\N$. By the formula \eqref{conc_nXde2}, for all $\e\in(0,1)$ and $d\in\N$ it is easy to obtain the following inequalities
\begin{eqnarray*}
\sum_{k=1}^{n^{X_d}(\e)-1}\lambda^{X_d}_k< (1-\e^2)\Lambda^{X_d} \leqslant \sum_{k=1}^{n^{X_d}(\e)}\lambda^{X_d}_k,
\end{eqnarray*}
which are reduced to the required bounds
\begin{eqnarray*}
n^{X_d}(\e)-1< (1-\e^2)\bigl(\Lambda^X/\lambda^X_1\bigr)^d \leqslant n^{X_d}(\e)\quad\text{for all}\quad \e\in(0,1),\,\, d\in\N,
\end{eqnarray*}
i.e. $n^{X_d}(\e)=\bigl\lceil\,(1-\e^2)\bigl(\Lambda^X/\lambda^X_1\bigr)^d\, \bigr\rceil$, where $\lceil\,\cdot\,\rceil$ is a ceiling function.

\subsection{Logarithmic asymptotics of the approximation complexity}\newparagraph
For tensor degree-type random elements stable distributions (see Appendix) play a crucial role in approximation problems, as the following theorem shows.

\begin{Theorem}\label{th_nXde_LimitDistr_TDRE}
Let $(a_d)_{d\in\N}$ be a sequence, $(b_d)_{d\in\N}$ be a positive  sequence such that $b_d\to\infty$, $d\to\infty$.  Let a non-increasing function $q\colon (0,1)\to\R$ and a distribution function $G$ satisfy the equation   $q(\e)=G^{-1}(1-\e^2)$ for all $\e\in\ContSet(q)$. Suppose that we have
\begin{eqnarray}\label{th_nXde_LimitDistr_TDRE_cond}
\ln n^{X_d}(\e)=a_d+q(\e)b_d+o(b_d),\quad d\to\infty,\quad\text{for all}\quad \e\in\ContSet(q).
\end{eqnarray}
Then $G$ is a stable distribution function. If $G$ is non-degenerate, then $b_{d+1}/b_d\to\infty$, $d\to\infty$. If $G$ is non-degenerate and non-normal, then $G=\StL_{\alpha,\rho,1,\mu}$  for some $\alpha\in(0,2)$, $\rho>0$, $\mu\in\R$.
\end{Theorem}
\textbf{Proof of Theorem \ref{th_nXde_LimitDistr_TDRE}.}\quad By Theorem \ref{th_nXde_LogAsymp_GRE}, the condition \eqref{th_nXde_LimitDistr_TDRE_cond} is equivalent to the convergence \eqref{th_nXde_LogAsymp_TPRE_conv}, where the functions $G^{X_d}_{a_d,b_d}$, $d\in\N$, are defined by \eqref{def_GXdadbd_GRE} for given  $a_d$, $b_d$, and $X_d$, $d\in\N$. According to Lemma \ref{lm_Uj}, $G^{X_d}_{a_d,b_d}(x)$ admits the representation \eqref{conc_GXdadbdSumUj_TPRE}, where the independent random variables $U_j$, $j\in\N$, have the common distribution:
\begin{eqnarray}\label{def_Uj_TDRE}
\Probab\bigl(U_j=\bigl|\ln \nlambda^{X}_k\bigr| \bigr)= \mult^{X} \bigl(\nlambda^{X}_k\bigr)\cdot\nlambda^{X}_k,\quad j\in\N,
\end{eqnarray}
where $k\in\N$ is such that $\lambda^{X}_k>0$. By Theorem \ref{th_LimitStable_Probab} (see Appendix),  the weak limit of $G^{X_d}_{a_d,b_d}(x)$ is necessarily stable.  

It is easy to check that the condition \eqref{cond_UN_TPRE} is satisfied for tensor degree-type random elements (set $X_{1,j}\colonequals X$) for any positive sequence $(b_d)_{d\in\N}$ such that $b_d\to\infty$, $d\to\infty$. Then Theorem \ref{th_nXde_LimitDistr_TPRE} gives the required assertion for $(b_d)_{d\in\N}$.

Assume that $G=\StL_{\alpha,\rho, \beta, \mu}$ for some $\alpha\in(0,2)$, $\rho>0$, $\beta\in[-1,1]$, and $\mu\in\R$. Let $L$ denote its L\'evy spectral function. From Theorem \ref{th_nXde_LimitDistr_TPRE} we have $L(x)=0$ for all $x<0$. In view of formula \eqref{conc_LevySpectrFuncStable_Probab}, we conclude that $\beta=1$.\quad $\Box$

According to the previous theorem, we can consider only functions $q$ such that $q(\e)=G^{-1}(1-\e^2)$ for all $\e\in\ContSet(q)$, where $G$ is a stable distribution function. Without loss of generality, we will restrict ourselves to the cases, where $G$ is non-degenerate and   it also has the standard form (see \eqref{conc_PhiStd_Probab} and \eqref{conc_StLStd_Probab}). 

Before formulating the criteria for \eqref{th_nXde_LimitDistr_TDRE_cond}, we consider conditions that appear there. The wide class of cases corresponds to the assumption:
\begin{eqnarray}\label{cond_2Mom_TDRE}
\sum_{k\in\N:\, \nlambda^X_k>0}\,\bigl|\ln \nlambda^X_k\bigr|^2\,\nlambda^X_k<\infty.
\end{eqnarray}
We also consider marginal random elements $X$ with the following regular variation of $\nlambda^X_k$, $k\in\N$:
\begin{eqnarray}\label{cond_RV_TDRE}
\sum_{k=1}^\infty\nlambda^X_k\,\id\bigl(\nlambda^X_k<e^{-x} \bigr)=x^{-\alpha} \varphi(x),
\end{eqnarray}
where $\alpha\geqslant0$ and $\svf$ is some slowly varying function at infinity (SVF for short, see Appendix). Under the assumption \eqref{cond_RV_TDRE} with some $\alpha>2$, the condition \eqref{cond_2Mom_TDRE} always holds. If either \eqref{cond_2Mom_TDRE} or \eqref{cond_RV_TDRE} with $\alpha>1$ hold, then the entropy of $\nlambda^X_k$, $k\in\N$, is well defined:
\begin{eqnarray}\label{def_EX_TDRE}
E^X:=\sum_{k\in\N:\, \nlambda^X_k>0}\bigl|\ln \nlambda^X_k\bigr|\,\nlambda^X_k.
\end{eqnarray} 
Under the assumption \eqref{cond_2Mom_TDRE}, the following deviation characteristic is also important:
\begin{eqnarray}\label{def_sigmaX_TDRE}
\sigma^X:=\biggl(\,\sum_{k\in\N:\, \nlambda^X_k>0}\,\Bigl(\bigl|\ln\nlambda^X_k\bigr|-E^X\Bigr)^2\,\nlambda^X_k\biggr)^{1/2}.
\end{eqnarray}

The next theorem provides a criterion of the asymptotics \eqref{th_nXde_LimitDistr_TDRE_cond}, where $q$ is a quantile of the distribution function $\Phi$ of the standard normal law.

\begin{Theorem}\label{th_nXde_LogAsymp_Normal_TDRE}
Let $(a_d)_{d\in\N}$ be a sequence, $(b_d)_{d\in\N}$ be a positive sequence such that $b_d\to\infty$, $d\to\infty$. For the asymptotics
\begin{eqnarray}\label{th_nXde_LogAsymp_Normal_TDRE_nXde}
\ln n^{X_d}(\e)=a_d+\Phi^{-1}(1-\e^2)b_d+o(b_d),\quad d\to\infty,\quad\text{for all}\quad\e\in(0,1),
\end{eqnarray}
it is necessary and sufficient to have:
\begin{itemize}
\item[$(i)$] the condition \eqref{cond_2Mom_TDRE} with $\sigma^X>0$ or the condition $\eqref{cond_RV_TDRE}$ with $\alpha=2$ and some SVF $\varphi;$
\item[$(ii)$] $a_d=E^X d+o(b_d)$, $d\to\infty;$
\item[$(iii)$] $\lim\limits_{d\to\infty} \dfrac{d}{b_d^2} \Bigl(M^X_2(b_d)- M^X_1(b_d)^2\Bigr)=1$,
\end{itemize}
where
\begin{eqnarray*}
M^X_p(x)\colonequals\sum_{k=1}^\infty\bigl|\ln\nlambda^X_k\bigr|^p\,\nlambda^X_k \id(\nlambda^X_k\geqslant e^{-x}),\quad x\geqslant0,\quad p\in\{1,2\}.
\end{eqnarray*}
If \eqref{cond_2Mom_TDRE} holds with $\sigma^X>0$ then $(iii)$ is equivalent to $b_d\sim \sigma^X d^{1/2}$, $d\to\infty$. If $\eqref{cond_RV_TDRE}$ holds with $\alpha=2$ and SVF $\varphi$ but  \eqref{cond_2Mom_TDRE} fails, then $(iii)$ is equivalent to $b_d\sim d^{1/2}\varphi_2(d) $, $d\to\infty$, with SVF $\svf_2$, which is defined by $\svf_2(d)\colonequals \sqrt{2}\,\bigl(1/\sqrt{\deHaansvf}\,\bigr)^\# (d^{1/2})$, $d\in\N$, where $(\,\cdot\,)^\#$ is the de Bruijn conjugation \footnote{see Appendix.} and SVF $\deHaansvf$ is defined by
\begin{eqnarray}\label{def_deHaansvf}
\deHaansvf(x)\colonequals \int\limits_0^x\dfrac{\svf(t)}{t} \dd t,\quad x\geqslant0.
\end{eqnarray}
\end{Theorem}
\textbf{Proof of Theorem \ref{th_nXde_LogAsymp_Normal_TDRE}.}\quad  The function $\Phi$ is absolutely continuous and strictly increasing on $\R$. Hence, by Theorem \ref{th_nXde_LogAsymp_GRE}, the condition \eqref{th_nXde_LogAsymp_Normal_TDRE_nXde} is equivalent to the convergence
\begin{eqnarray}\label{conc_GXdadbdtoPhi_TDRE}
\lim_{d\to\infty} G^{X_d}_{a_d,b_d}(x)= \Phi(x)\quad\text{for all}\quad x\in \R,
\end{eqnarray}
where the functions $G^{X_d}_{a_d,b_d}$, $d\in\N$, are defined by \eqref{def_GXdadbd_GRE} for given  $a_d$, $b_d$, and $X_d$, $d\in\N$.  By \eqref{conc_GXdadbdSumUj_TPRE}, $G^{X_d}_{a_d,b_d}(x)$, $d\in\N$, are distribution functions of centered and normalized sums of independent random variables $U_j$, $j\in\N$, with the common distribution \eqref{def_Uj_TDRE}.

Let us show the sufficiency of the assumptions $(i)$--$(iii)$ for the convergence \eqref{conc_GXdadbdtoPhi_TDRE}. In probabilistic interpretation the assumption \eqref{cond_2Mom_TDRE} means $\Expec U_j^2<\infty$. Also we can rewrite the condition \eqref{cond_RV_TDRE} for $\alpha=2$ as follows: $\Probab(U_j>x)= x^{-2}\svf(x)$. Since $U_j\geqslant 0$, we have $\Probab(U_j<-x)=0$ for all $x>0$. Then the conditions $(i)$ and $(iii)$ are sufficient for the convergence $\lim\limits_{d\to\infty}G^{X_d}_{E^X d, b_d}(x)=\Phi(x)$, $x\in\R$, by Theorem \ref{th_ConvtoNormal_Probab}. Therefore, using the assumption $(ii)$ of the theorem and the assertion $(ii)$ of Lemma \ref{lm_anbnWeakConv_Probab} we obtain \eqref{conc_GXdadbdtoPhi_TDRE}.

Let us show the necessity of $(i)$--$(iii)$. Under the convergence \eqref{conc_GXdadbdtoPhi_TDRE}, the condition $(i)$ holds by Theorem \ref{th_ConvtoNormal_Probab}. Also we have $\lim\limits_{d\to\infty}G^{X_d}_{E^X d, b^*_d}(x)=\Phi(x)$ for any $x\in\R$ and some sequence $(b^*_d)_{d\in\N}$ that satisfies $(iii)$. According to the assertion $(i)$ of Lemma \ref{lm_anbnWeakConv_Probab}, we obtain $b^*_d\sim b_d$ and $a_d=E^X d +o(b^*_d)$, $d\to\infty$, that directly yields $(ii)$. The condition $(iii)$ for $(b_d)_{d\in\N}$ follows from the slow variation of the function $x\mapsto M^X_2(x)-M^X_1(x)^2$, $x\geqslant0$ (which is  justified the next notes). 

Under the condition $\eqref{cond_2Mom_TDRE}$, we have $M^X_2(x)-M^X_1(x)^2\to\sigma^X$, $x\to\infty$, i.e. $(iii)$ is equivalent to $b_d\sim\sigma^X d^{1/2}$, $d\to\infty$. In the remainder of this proof we assume  $\eqref{cond_RV_TDRE}$ with $\alpha=2$ and some SVF $\svf$ but \eqref{cond_2Mom_TDRE} fails. Here we have $M^X_2(x)\to\infty$, $x\to\infty$, and, by the remark (2.6.14) from \cite{IbrLin} (see p. 80),  $M^X_1(x)^2=o\bigl(M^X_2(x)\bigr)$, $x\to\infty$. Hence $(iii)$ is equivalent to $b_d^2\sim d M^X_2(b_d)$, $d\to\infty$. Using integral representation for $M^X_2(x)$ and integrating it by parts, we get
\begin{eqnarray*}
M^X_2(x)=\int\limits_0^{x} t^2\dd\Bigl(\sum_{k=1}^\infty\nlambda^X_k \id(\nlambda^X_k\geqslant e^{-t}) \Bigr)=\int\limits_0^{x} t^2\dd \bigl(1-t^{-2}\svf(t)\bigr)=-\svf(x)+2 \deHaansvf(x),\quad x\geqslant 0,
\end{eqnarray*}
where $\deHaansvf$ is defined by \eqref{def_deHaansvf}. By Lemma \ref{lm_deHaan_SVF} (see Appendix), $\deHaansvf$ is a SVF and also $\svf(x)=o(\deHaansvf(x))$ as $x\to\infty$. Hence $(iii)$ is equivalent to $b_d^2\sim 2 d\,\deHaansvf(b_d)$, $d\to\infty$. Rewriting this in the equivalent form
\begin{eqnarray*}
d\sim \bigl(b_d/\sqrt{2}\,\bigr)^2\cdot \bigl(1/\sqrt{\deHaansvf}\,\bigr)^2\bigl(b_d/\sqrt{2}\, \bigr),\quad d\to\infty,
\end{eqnarray*}
and using Lemma \ref{lm_AsymInv_deBruijn_SVF} from Appendix, we obtain the required relation for $b_d$, $d\in\N$.\quad $\Box$\\

Using tools of regular variation theory (see \cite{BingGoldTeug}), it is possible to find simpler asymptotic versions of the function $\svf_2$ from Theorem \ref{th_nXde_LogAsymp_Normal_TDRE} under special assumptions (see Lemma \ref{lm_Bekessy_SVF} and examples in the next subsection).

The previous theorem has the following important corollary concerning wide class of tensor degree-type random elements.
\begin{Theorem}\label{th_nXde_LogAsymp_Normal_2Mom_TDRE}
Under the assumption \eqref{cond_2Mom_TDRE} with $\sigma^X>0$, we have
\begin{eqnarray*}	
\ln n^{X_d}(\e)=E^X d+\Phi^{-1}(1-\e^2)\sigma^X d^{1/2}+o(d^{1/2}),\quad d\to\infty,\quad\text{for all}\quad\e\in(0,1).
\end{eqnarray*}
\end{Theorem}
In fact, this theorem was obtained by M. A. Lifshits and E. V. Tulyakova in the paper \cite{LifTul}. However, strictly speaking, the proof from \cite{LifTul} was done only for sequences $(\nlambda^X_k)_{k\in\N}$ with unit multiplicity of every element, i.e. $\mult^X(\nlambda^X_k)=1$, $k\in\N$ (it is hidden in the last formula on p. 106 in \cite{LifTul}).

Next remarks  follow directly from the definitions of $E^X$ and $\sigma^X$.
\begin{Remark}
The condition \eqref{cond_lambdaXk_triv_TDRE} holds iff $\sigma^X=0$.
\end{Remark}
\begin{Remark}
The equality $\lambda^X_1=\Lambda^X$ $($i.e. \eqref{cond_lambdaXk_triv_TDRE} with $l^X=1$$)$ holds iff $E^X=0$.
\end{Remark}
These remarks show that there is no loss of generality for us in assuming $\sigma^X>0$ in Theorems \ref{th_nXde_LogAsymp_Normal_TDRE} and~\ref{th_nXde_LogAsymp_Normal_2Mom_TDRE}. Also we can now conclude that under the assumption $(i)$ of the previous theorem, the complexty $n^{X_d}(\e)$ grows mainly \textit{exponentially} with the constant $E^X>0$ as $d\to\infty$.

The next theorem provides a criterion of the asymptotics \eqref{th_nXde_LimitDistr_TDRE_cond} with $G=\StL_{\alpha,1}$, $\alpha\in(0,2)$.
\begin{Theorem}\label{th_nXde_LogAsymp_Stable_TDRE}
Let $(a_d)_{d\in\N}$ be a sequence, $(b_d)_{d\in\N}$ be a positive sequence such that $b_d\to\infty$, $d\to\infty$. For the asymptotics
\begin{eqnarray}\label{th_nXde_LogAsymp_Stable_TDRE_nXde}
\ln n^{X_d}(\e)=a_d+\StL^{-1}_{\alpha,1}(1-\e^2)b_d+o(b_d),\quad d\to\infty,\quad\text{for all}\quad\e\in(0,1).
\end{eqnarray}
It is necessary and sufficient to have:
\begin{itemize}
\item[$(i)$]the condition $\eqref{cond_RV_TDRE}$ with given $\alpha$ and some SVF $\varphi;$
\item[$(ii)$] $a_d=a^*_d+o(b_d)$, $d\to\infty;$
\item[$(iii)$] $\lim\limits_{d\to\infty}d  \sum\limits_{k=1}^\infty\nlambda^X_k\,\id\bigl(\nlambda^X_k<e^{-b_d}\bigr)=1$.
\end{itemize}
Here $a^*_d\colonequals0$ for $\alpha\in (0,1)$, $a^*_d\colonequals d E^X$ for $\alpha\in (1,2)$ and 
\begin{eqnarray*}
a^*_d\colonequals  d \sum_{k=1}^{\infty} |\ln\nlambda^X_k|\,\nlambda^X_k\id\bigl(\nlambda^X_k>e^{-b_d}\bigr) +(1-\EulerCon) b_d,\quad \text{for}\quad \alpha=1,
\end{eqnarray*}
where $\EulerCon$ is the Euler constant. Under the assumption $(i)$ with some $\alpha\in(0,2)$, the condition $(iii)$ is equivalent to $b_d\sim d^{1/\alpha} \svf_\alpha(d)$, $d\to\infty$, with SVF $\svf_\alpha$, defined by $\svf_\alpha(d)\colonequals\bigl((1/\svf)^{1/\alpha}\bigr)^\#(d^{1/\alpha})$, $d\in\N$, where $(\,\cdot\,)^\#$ is the de Bruijn conjugation.
\end{Theorem}
\textbf{Proof of Theorem \ref{th_nXde_LogAsymp_Stable_TDRE}.}\quad Since $\StL_{\alpha,1}$ is self-decomposable (see Appendix), the distribution function $\StL_{\alpha,1}$ is absolutly continuous on $\R$ and strictly increasing on $(\lext \StL_{\alpha,1}, \rext \StL_{\alpha,1})$ in view of Remark \ref{rem_StrIncSelfDec_Probab}. Hence, by Theorem \ref{th_nXde_LogAsymp_GRE}, the condition \eqref{th_nXde_LogAsymp_Stable_TDRE_nXde} is equivalent to the convergence
\begin{eqnarray}\label{conc_GXdadbdtoStL_TDRE}
\lim_{d\to\infty} G^{X_d}_{a_d,b_d}(x)= \StL_{\alpha,1}(x)\quad\text{for all}\quad x\in \R,
\end{eqnarray}
where the functions $G^{X_d}_{a_d,b_d}$, $d\in\N$, are defined by \eqref{def_GXdadbd_GRE} for given  $a_d$, $b_d$, and $X_d$, $d\in\N$.  By probabilistic representation  \eqref{conc_GXdadbdSumUj_TPRE}, $G^{X_d}_{a_d,b_d}(x)$, $d\in\N$, are distribution functions of centered and normalized sums of independent random variables $U_j$, $j\in\N$, with the common distribution \eqref{def_Uj_TDRE}.  

Let us show the sufficiency of the assumptions $(i)$--$(iii)$ for the convergence \eqref{conc_GXdadbdtoStL_TDRE}. In probabilistic interpretation the assumption \eqref{cond_RV_TDRE} can be written as follows: $\Probab(U_j>x)= x^{-\alpha}\svf(x)$, $x\geqslant 0$. Since $U_j\geqslant 0$, we have $\Probab(U_j<-x)=0$ for all $x>0$. By Theorem \ref{th_ConvtoStable_Probab}, the conditions $(i)$ and $(iii)$ are sufficient for the convergence $\lim\limits_{d\to\infty}G^{X_d}_{a^*_d, b_d}(x)=\StL_{\alpha,1}(x)$, $x\in\R$. Therefore, using the assumption $(ii)$ of the theorem and the assertion $(ii)$ of Lemma \ref{lm_anbnWeakConv_Probab} we obtain \eqref{conc_GXdadbdtoStL_TDRE}.

Let us show the necessity of $(i)$--$(iii)$. Under the convergence \eqref{conc_GXdadbdtoStL_TDRE}, the condition $(i)$ holds by Theorem \ref{th_ConvtoStable_Probab}. Also we have $\lim\limits_{d\to\infty}G^{X_d}_{a^*_d, b^*_d}(x)=\StL_{\alpha,1}(x)$ for any $x\in\R$ and some sequence $(b^*_d)_{d\in\N}$ that satisfies $(iii)$. According to the assertion $(i)$ of Lemma \ref{lm_anbnWeakConv_Probab}, we obtain $b^*_d\sim b_d$ and $a_d=a^*_d +o(b^*_d)$, $d\to\infty$, that  yields $(ii)$. The condition $(iii)$ for $(b_d)_{d\in\N}$ follows from the regular variation of the function $x\mapsto \sum\limits_{k=1}^\infty\lambda^X_k\,\id\bigl(\lambda^X_k<e^{-x}\bigr)$, $x\geqslant0$. 

Under the condition $\eqref{cond_RV_TDRE}$ with some $\alpha\in(0,2)$ and some SVF $\svf$, $(iii)$ can be rewritten as $d\sim (b_d)^\alpha\,(1/\svf)^{1/\alpha}(b_d)^{\alpha}$, $d\to\infty$. Using Lemma \ref{lm_AsymInv_deBruijn_SVF} from Appendix, we obtain the required relation for $b_d$, $d\in\N$.\quad $\Box$\\

Simpler asymptotic versions of $b_d$ can be obtained using Lemma \ref{lm_Bekessy_SVF} from Appendix.

\begin{Remark}
Under the assumption $\eqref{cond_RV_TDRE}$ with $\alpha=1$, $(ii)$, and $(iii)$, we have $b_d=o(a_d)$, $d\to\infty$.
\end{Remark}
Indeed, using integral representation for $a_d$ and the condition $(iii)$, we obtain as $d\to\infty$
\begin{eqnarray*}
a_d&=&d \int\limits_{0}^{b_d} t \dd \Bigl(\sum\limits_{k=1}^\infty\lambda^X_k\,\id\bigl(\lambda^X_k\geqslant e^{-t}\bigr) \Bigr)+ (1-\EulerCon)b_d+o(b_d)\\
&=&d \int\limits_{0}^{b_d} t \dd (1-t^{-1}\svf(t))+ (1-\EulerCon)b_d+o(b_d)\\
&=&-d\svf (b_d)+d\deHaansvf(b_d)+(1-\EulerCon)b_d+o(b_d)\\
&=& d\deHaansvf(b_d)-\EulerCon b_d+o(b_d)\\
&=&b_d\bigl(\deHaansvf(b_d)/\svf(b_d)\bigr)-\EulerCon b_d+o(b_d).
\end{eqnarray*}
By Lemma \ref{lm_deHaan_SVF} (see Appendix), we have $\svf(x)=o(\deHaansvf(x))$ as $x\to\infty$, which  completes the verification.

Let us analyze the behaviour of $n^{X_d}(\e)$ as $d\to\infty$. Under the assumption $(i)$ with $\alpha\in(1,2)$, the approximation complexity grows \textit{exponentially} due to $a_d$ with \textit{subexponential} term containing $b_d$. In the boundary case $\alpha=1$ the main term $a_d$ can yield more than exponential growth (see example in the next subsection). Under the cases $\alpha\in(0,1)$, factor $d^{1/\alpha}$ (of $b_d$) gives the \textit{superexponential} growth of the quantity $n^{X_d}(\e)$.

Consider the unexplored pathological case when the assumption \eqref{cond_RV_TDRE} is satisfied for $\alpha=0$ and some SVF $\svf$ that $\svf(x)\to0$ as $x\to\infty$. Here it is impossible to find $(a_d)_{d\in\N}$ and $(b_d)_{d\in\N}$ for obtaining the asymptotics \eqref{th_nXde_LimitDistr_TDRE_cond} with non-degenerate distribution function $G$ (see comments in Subsection \ref{subsec_Darling} in Appendix). Nevertheless, we obtain the following result.

\begin{Theorem}\label{th_nXde_LogAsymp_SV_TDRE}
If $\eqref{cond_RV_TDRE}$ holds with $\alpha=0$ and some SVF $\svf$, then
\begin{eqnarray}\label{th_nXde_LogAsymp_SV_TDRE_conc}
\lim_{d\to\infty}d\,\svf\bigl(\ln n^{X_d}(\e)\bigr)=-\ln(1-\e^2)\quad\text{for all}\quad\e\in(0,1).
\end{eqnarray}
\end{Theorem}
\textbf{Proof of Theorem \ref{th_nXde_LogAsymp_SV_TDRE}.}\quad At first we will prove that
\begin{eqnarray}\label{conc_nlambdaavg_SV_TDRE}
\lim_{d\to\infty}\Bigl[d\,\svf\Bigl(\bigl|\ln\nlambda^{X_d}(\e)\bigr|\Bigr)\Bigr]=-\ln(1-\e^2)\quad\text{for all}\quad\e\in(0,1),
\end{eqnarray}
where $\nlambda^{X_d}(\e)$ is defined by \eqref{def_nlambdaXde}. Let $U_j$, $j\in\N$, be independent random variables with common distribution \eqref{def_Uj_TDRE}. Then we have
\begin{eqnarray*}
\Probab(U_1>x)=\sum_{k=1}^\infty\nlambda^X_k\,\id\bigl(\nlambda^X_k<e^{-x}\bigr)=\svf(x), \quad x\geqslant0.
\end{eqnarray*}
Introduce the function $\Laplsvf$:
\begin{eqnarray*}
\Laplsvf(x)\colonequals 1-\Expec \exp\{-U_1/x\}= 1-\sum_{k=1}^\infty\bigl(\nlambda^X_k\bigr)^{1+1/x} , \quad x>0.
\end{eqnarray*}
By remarks to Theorem \ref{th_Darling2_Probab} (see Appendix), the function $\Laplsvf$ is continuous and strictly decreasing on $(0,\infty)$. Also it satisfies:
\begin{eqnarray}\label{conc_Laplsvf_equiv_SV_TDRE}
\Laplsvf(x)\sim\svf(x), \quad x\to\infty.
\end{eqnarray}

Let us set
\begin{eqnarray*}
F^{X_d}(x)\colonequals\sum_{k=1}^\infty\nlambda^{X_d}_k\,\id\Bigl(-d\,\Laplsvf\bigl(\bigl|\ln\nlambda^{X_d}_k\bigr|\bigr) \leqslant x\Bigr), \quad x\in \R.
\end{eqnarray*}
According to Lemma \ref{lm_Uj}, $F^{X_d}$, $d\in\N$, have the following representations:
\begin{eqnarray*}
F^{X_d}(x)=\Probab\biggl(-d\,\Laplsvf\Bigl(\sum_{j=1}^d U_j\Bigr)\leqslant x\biggr), \quad x\in \R.
\end{eqnarray*}
By Theorem \ref{th_Darling2_Probab}, for $x<0$ we get:
\begin{eqnarray*}
	F^{X_d}(x)=1-\Probab\biggl(d\,\Laplsvf\Bigl(\sum_{j=1}^d U_j\Bigr)< -x\biggr)\to e^{x}, \quad d\to\infty.
\end{eqnarray*}
From this we infer the convergence $(F^{X_d})^{-1}(r)\to \ln r$, $d\to\infty$, for any $r\in(0,1)$. In particular, we have
\begin{eqnarray*}
\lim_{d\to\infty}(F^{X_d})^{-1}(1-\e^2)= \ln(1-\e^2)\quad\text{for all}\quad\e\in(0,1).
\end{eqnarray*}
It is true that $(F^{X_d})^{-1}(1-\e^2)=-d\,\Laplsvf\bigl(\bigl|\ln\nlambda^{X_d}(\e)\bigr|\bigr)$ for any $\e\in(0,1)$ and $d\in\N$. Indeed, by \eqref{conc_nlambdaXde_inf}, we have
\begin{eqnarray*}
	d\,\Laplsvf\bigl(\bigl|\ln\nlambda^{X_d}(\e)\bigr|\bigr)&=&d\,\Laplsvf\Bigl(\inf\bigl\{x\in\R: \sum_{k=1}^\infty\nlambda^{X_d}_k\,\id\bigl(\bigl|\ln\nlambda^{X_d}_k\bigr| \leqslant x\bigr)\geqslant 1-\e^2\bigr\}\Bigr)\\
	&=&d\,\Laplsvf\biggl(\inf\Bigl\{x\in\R: \sum_{k=1}^\infty \nlambda^{X_d}_k\,\id\Bigl(d\,\Laplsvf\bigl(\bigl|\ln\nlambda^{X_d}_k\bigr|\bigr) \geqslant d\,\Laplsvf(x)\Bigr)\geqslant 1-\e^2\Bigr\}\biggr)\\
	&=&-\inf\bigl\{y\in\R: F^{X_d}(y)\geqslant 1-\e^2\bigr\}\\
	&=&-(F^{X_d})^{-1}(1-\e^2).
\end{eqnarray*}
Hence for any $\e\in(0,1)$ we obtain $\lim\limits_{d\to\infty} d\,\Laplsvf\bigl(\bigl|\ln\nlambda^{X_d}(\e)\bigr|\bigr)=-\ln(1-\e^2)$. By the strict decay of $\Laplsvf$, for any $\e\in(0,1)$  $\bigl|\ln\nlambda^{X_d}(\e)\bigr|\to \infty$ as $d\to\infty$. 
The convergence \eqref{conc_nlambdaavg_SV_TDRE} follows from \eqref{conc_Laplsvf_equiv_SV_TDRE}.

Next, we prove \eqref{th_nXde_LogAsymp_SV_TDRE_conc}. From \eqref{conc_nXde_nlambdaXde_upbound} and strict decay of $\Laplsvf$ it follows that
\begin{eqnarray*}
	\varliminf_{d\to\infty} d\,\Laplsvf\bigl(\ln n^{X_d}(\e)\bigr)\geqslant\lim_{d\to\infty} d\,\Laplsvf\bigl(\bigl|\ln\nlambda^{X_d}(\e)\bigr|\bigr)=-\ln(1-\e^2).
\end{eqnarray*}
Fix arbitrary $h>0$ and $c_h\in(1,1/\e)$ such that 
$\ln(1-c_h^2\e^2)/\ln(1-\e^2)<e^h$. Using \eqref{conc_nXde_nlambdaXde_lowbound} with $\e_1=\e$ and $\e_2=c_h\e$, $\bigl|\ln\nlambda^{X_d}(\e)\bigr|\to \infty$, $d\to\infty$, and slow variation of $\Laplsvf$, we obtain:
\begin{eqnarray*}
\varlimsup_{d\to\infty} d\,\Laplsvf\bigl(\ln n^{X_d}(\e)\bigr)&\leqslant& \lim_{d\to\infty} d\,\Laplsvf\bigl(\bigl|\ln\nlambda^{X_d}(\e)\bigr|  + \ln\bigl((c_h^2-1)\e^2/2\bigr)\bigr)\\
&=& \lim_{d\to\infty} d\,\Laplsvf\bigl(\bigl|\ln\nlambda^{X_d}(\e)\bigr| \bigr)\\
&=&-\ln(1-c_h^2\e^2)\\
&\leqslant& -\ln(1-\e^2) e^h.
\end{eqnarray*}
Hence  $\lim\limits_{d\to\infty} d\,\Laplsvf\bigl(\ln n^{X_d}(\e)\bigr)=-\ln(1-\e^2)$ for any $\e\in(0,1)$. The equivalence \eqref{conc_Laplsvf_equiv_SV_TDRE} yields \eqref{th_nXde_LogAsymp_SV_TDRE_conc}. \quad $\Box$

\begin{Corollary}\label{co_th_nXde_LogAsymp_SV}
For any $\e\in(0,1)$  $\bigl(\ln n^{X_d}(\e)\bigr)_{d\in\N}$ is rapidly varying sequence, i. e.
\begin{eqnarray*}
\lim\limits_{d\to\infty} \dfrac{\ln n^{X_{ \lfloor cd\rfloor}}(\e)}{\ln n^{X_d}(\e)}=\infty\quad\text{for all}\quad c>1,
\end{eqnarray*}
where $\lfloor\cdot\rfloor$ is a floor function.
\end{Corollary}
Thus here growth of $n^{X_d}(\e)$  is extremely fast as $d\to\infty$ (it can be double exponential, $\exp\{\exp\{\,\cdot\,\}\}$, see an example in the next subsection).

To prove this corollary we suppose, contrary to our claim that for some $c>1$ there exists a subsequence $(d_k)_{k\in\N}$ such that
\begin{eqnarray}\label{conc_assumlamavgcde_SV_TDRE}
\lim_{k\to\infty}\dfrac{\ln n^{X_{\lfloor c d_k\rfloor}}(\e)}{\ln n^{X_{d_k}}(\e)}=p>0.
\end{eqnarray} 
On the one hand,  by Theorem \ref{th_nXde_LogAsymp_SV_TDRE}, we obtain
\begin{eqnarray*}
	\lim_{k\to\infty}\Bigl(\lfloor c d_k\rfloor\,\svf\bigl(\ln n^{X_{\lfloor c d_k\rfloor}}(\e))\bigr)\Bigr)=-\ln(1-\e^2).
\end{eqnarray*}
On the other hand, from slow variation of $\svf$ and \eqref{conc_assumlamavgcde_SV_TDRE} it follows
\begin{eqnarray*}
	\lim_{k\to\infty}\lfloor c d_k\rfloor\,\svf\bigl(\ln n^{X_{\lfloor c d_k\rfloor}}(\e))\bigr)=\lim_{k\to\infty}\lfloor c d_k\rfloor\,\svf\bigl(\ln n^{X_{d_k}}(\e)\bigr)=-c\ln(1-\e^2).
\end{eqnarray*}
This leads to a contradiction.

\subsection{Applications}\newparagraph
Let us consider a sequence of tensor degrees $X_d=X^{\otimes d}$, $d\in\N$. Suppose that eigenvalues of $K^X$ has the following asymptotics:
\begin{eqnarray}\label{cond_lambdaXk_RegBehav}
\dfrac{\lambda^X_k}{\Lambda^X}\sim \dfrac{\beta}{k^p (\ln k)^{1+r}},\quad k\to\infty,
\end{eqnarray}
where $\beta>0$ and numbers $p$ and $r$ must satisfy $\sum_{k=1}^{\infty} \lambda^X_k<\infty$, i.e. $p>1$, $r\in\R$ or $p=1$, $r>0$.

The following  assertion is a direct corollary of Theorem \ref{th_nXde_LogAsymp_Normal_2Mom_TDRE}.
\begin{Proposition}
Under the assumption \eqref{cond_lambdaXk_RegBehav} with $\beta>0$, $p>1$, $r\in\R$, or $p=1$, $r>2$, we have
\begin{eqnarray*}
\ln n^{X_d}(\e)=E^X d+\Phi^{-1}(1-\e^2)\sigma^X d^{1/2}+o(d^{1/2}),\quad d\to\infty,\quad\text{for all}\quad \e\in(0,1),
\end{eqnarray*}
where $E^X$ and $\sigma^X$ are defined by \eqref{def_EX_TDRE} and \eqref{def_sigmaX_TDRE}, respectively.
\end{Proposition}

In order to consider the remaining cases  we need the following auxiliar lemma.
\begin{Lemma}\label{lm_sumnlambdaXk_RegBehav}
Under the assumption \eqref{cond_lambdaXk_RegBehav} with $p=1$, $r>0$, we have
\begin{eqnarray*}
\sum_{k=1}^{\infty} \nlambda^X_k \id\bigl(\nlambda^X_k<e^{-x}\bigr)\sim \dfrac{\beta}{r}\cdot x^{-r},\quad x\to\infty. 
\end{eqnarray*}
\end{Lemma}
\textbf{Proof of Lemma \ref{lm_sumnlambdaXk_RegBehav}.}\quad First, we get
\begin{eqnarray}\label{conc_sumnlambdaXkn_RegBehav}
\sum_{k=n}^{\infty} \nlambda^X_k\sim\sum_{k=n}^{\infty} \dfrac{\beta}{k\, (\ln k)^{1+r}}\sim \int\limits_{n}^{\infty} \dfrac{\beta\dd t}{t (\ln t)^{1+r}} = \dfrac{\beta}{r}\cdot(\ln n)^{-r}, \quad n\to\infty.
\end{eqnarray}
Set $k_x\colonequals \min\bigl\{k\in\N: \nlambda^X_k<e^{-x}\bigr\}$. Using Lemmas \ref{lm_AsymInv_deBruijn_SVF}  and \ref{lm_Bekessy_SVF} (see Appendix), we find
\begin{eqnarray}\label{conc_kx_asymp_RegBehav}
k_x\sim e^x \bigl((\ln)^{1+r}/\beta\bigr)^\#(e^x)\sim \dfrac{\beta e^x}{x^{1+r}},\quad x\to\infty,
\end{eqnarray}
and, in particular, $\ln k_x\sim x$, $x\to\infty$. From this and \eqref{conc_sumnlambdaXkn_RegBehav}
we obtain the required asymptotics.\quad $\Box$

\begin{Proposition}\label{pr_nXde_p1r2_RegBehav}
Under the assumption \eqref{cond_lambdaXk_RegBehav} with $\beta>0$, $p=1$, and $r=2$, for all $\e\in(0,1)$ we have 
\begin{eqnarray*}
\ln n^{X_d}(\e)=E^X d+\Phi^{-1}(1-\e^2)\,(\beta/2)^{1/2} (d \ln d)^{1/2}+o((d \ln d)^{1/2}),\quad d\to\infty,
\end{eqnarray*}
where $E^X$ is defined by \eqref{def_EX_TDRE} .
\end{Proposition}
\textbf{Proof of Proposition \ref{pr_nXde_p1r2_RegBehav}.}\quad On account of Lemma  \ref{lm_sumnlambdaXk_RegBehav} for $r=2$, the required asymptotics is obtained by Theorem \ref{th_nXde_LogAsymp_Normal_TDRE} with $\svf(x)\sim \beta/2$, $x\to\infty$, and $a_d=E^X d$, $b_d\colonequals d^{1/2}\sqrt{2}\,\bigl(1/\sqrt{\deHaansvf}\,\bigr)^\# (d^{1/2})$, where $\deHaansvf$ is defined by \eqref{def_deHaansvf}. We only need to find asymptotics of the sequence $(b_d)_{d\in\N}$. It easily seen that
\begin{eqnarray*}
\deHaansvf(x)= \int\limits_{0}^{x} \dfrac{\deHaansvf(t)}{t}\dd t\sim \int\limits_{1}^{x} \dfrac{\beta}{2 t}\dd t\sim \dfrac{\beta\ln x}{2},\quad x\to\infty.
\end{eqnarray*}
Hence, from Lemma \ref{lm_Bekessy_SVF} (see Appendix) we have
\begin{eqnarray*}
\bigl(1/\sqrt{\deHaansvf}\,\bigr)^\# (d^{1/2})\sim \deHaansvf(d^{1/2})^{1/2}\sim ((\beta/4)\ln d)^{1/2},\quad d\to\infty.
\end{eqnarray*}
This yields $b_d\sim (\beta/2)^{1/2} (d \ln d)^{1/2}$, $d\to\infty$. \quad $\Box$

\begin{Proposition}\label{pr_nXde_p1r02_RegBehav}
Under the assumption \eqref{cond_lambdaXk_RegBehav} with $\beta>0$, $p=1$, and $r\in(0,2)$, we have
\begin{eqnarray*}
\ln n^{X_d}(\e)=a_d+\StL_{r,1}^{-1}(1-\e^2)\,(\beta/r)^{1/r} d^{1/r}+o(d^{1/r}),\quad d\to\infty,\quad\text{for all}\quad\e\in(0,1),
\end{eqnarray*}
where $a'_d\colonequals0$ for $r\in (0,1)$, $a'_d\colonequals d E^X$ for $r\in (1,2)$, and 
\begin{eqnarray}\label{pr_nXde_p1r02_RegBehav_ad_r1}
a_d\colonequals  d \sum_{k=1}^{\infty} |\ln\nlambda^X_k|\,\nlambda^X_k\id\bigl(\nlambda^X_k>e^{-\beta d}\bigr) +(1-\EulerCon) \beta d,\quad \text{for}\quad r=1.
\end{eqnarray}
Here $\EulerCon$ is the Euler constant and $E^X$ is defined by \eqref{def_EX_TDRE}. In the case $r=1$  $a_d\sim \beta d\ln\ln d$, $d\to\infty$.
\end{Proposition}
\textbf{Proof of Proposition \ref{pr_nXde_p1r02_RegBehav}.}\quad According to Lemma \ref{lm_sumnlambdaXk_RegBehav} with given $r\in(0,2)$, we can use Theorem \ref{th_nXde_LogAsymp_Normal_TDRE} with $\alpha=r$ and $\svf(x)\sim \beta/r$, $x\to\infty$. 
Since, by Lemma \ref{lm_Bekessy_SVF}, we have 
\begin{eqnarray*}
	d^{1/r} \bigl((1/\svf)^{1/r}\bigr)^\#(d^{1/r})\sim d^{1/r} \svf(d^{1/r})^{1/r}\sim (\beta/r)^{1/r} d^{1/r},\quad d\to\infty,
\end{eqnarray*}
on account of  Lemma \ref{lm_anbnWeakConv_Probab}, numbers $b_d$ can be choosen as follows $b_d\colonequals(\beta/r)^{1/r} d^{1/r}$, $d\in\N$. The expressions for $a_d$ are directly obtained from Theorem \ref{th_nXde_LogAsymp_Normal_TDRE}. In particular, for the case $r=1$ we have the formula \eqref{pr_nXde_p1r02_RegBehav_ad_r1}. On account of \eqref{conc_kx_asymp_RegBehav}, we get $k_{\beta d} \colonequals\min\bigl\{k\in\N: \nlambda^X_k< e^{-\beta d} \bigr\}\sim e^{\beta d}/(\beta d^2)$. Accordingly, using assumptions on $\nlambda^X_k$, $k\in\N$,  it follows that 
\begin{eqnarray*}
	\sum_{k=1}^{\infty} |\ln\nlambda^X_k|\,\nlambda^X_k\id\bigl(\nlambda^X_k\geqslant e^{-\beta d}\bigr)\sim \sum_{k=1}^{k_{\beta d}} \dfrac{\beta}{k \ln k}\sim \ln\ln k_{\beta d} \sim\beta \ln\ln d, \quad d\to\infty. 
\end{eqnarray*}
This gives $a_d\sim \beta d \ln\ln d$, $d\to\infty$, for the case $r=1$.\quad $\Box$

In order to have full expansion of $a_d$ (up to $o(d)$) in Proposition \ref{pr_nXde_p1r02_RegBehav} for the case $r=1$, we need know more information about asymptotic behaviour of  $\lambda^X_k$, $k\in\N$. We omit these detailes, which reduce to routine calculations.

Now we consider an example of applying Theorem \ref{th_nXde_LogAsymp_SV_TDRE}.

\begin{Proposition}\label{pr_nXde_s_RegBehav}
Suppose that
\begin{eqnarray*}
\dfrac{\lambda^X_k}{\Lambda^X}\sim \dfrac{\beta}{k (\ln k) (\ln\ln k)^{1+s}},\quad k\to\infty,
\end{eqnarray*}
with $\beta>0$ and $s>0$. Then
\begin{eqnarray*}
\ln\ln n^{X_d}(\e)\sim \biggl(\dfrac{\beta d}{s |\ln(1-\e^2)|} \biggr)^{1/s},\quad d\to\infty,\quad\text{for all}\quad \e\in(0,1).
\end{eqnarray*}
\end{Proposition}
\textbf{Proof of Proposition \ref{pr_nXde_s_RegBehav}.}\quad First, we get
\begin{eqnarray*}
	\sum_{k=n}^{\infty} \nlambda^X_k\sim\sum_{k=n}^{\infty} \dfrac{\beta}{k(\ln k) (\ln\ln k)^{1+s}}\sim \int\limits_{n}^{\infty} \dfrac{\beta\dd t}{t\ln t (\ln\ln t)^{1+s}} = \dfrac{\beta}{s}\cdot(\ln\ln n)^{-s}, \quad n\to\infty.
\end{eqnarray*}
Set $k_x\colonequals \min\bigl\{k\in\N: \nlambda^X_k<e^{-x}\bigr\}$. Using Lemma \ref{lm_AsymInv_deBruijn_SVF} (see Appendix) and Lemma \ref{lm_Bekessy_SVF}, we find
\begin{eqnarray*}
	k_x\sim e^x \bigl((\ln )(\ln\ln)^{1+s}/\beta\bigr)^\#(e^x)\sim \dfrac{\beta e^x}{x (\ln x)^{1+s}},\quad x\to\infty,
\end{eqnarray*}
that gives $\ln\ln k_x\sim \ln x$, $x\to\infty$. From this we have
\begin{eqnarray*}
	\sum_{k=1}^{\infty} \nlambda^X_k\,\id\bigl(\nlambda^X_k<e^{-x}\bigr)= \sum_{k=k_x}^{\infty} \nlambda^X_k\sim\dfrac{\beta}{s}\cdot(\ln x)^{-s}, \quad x\to\infty.
\end{eqnarray*}
The required asymptotics conclude from Theorem \ref{th_nXde_LogAsymp_SV_TDRE}.\quad $\Box$

\section{Appendix}\newparagraph
Here we recall the definitions and basic properties of self-decomposable and stable probability distributions. Also we provide known limit theorems of weak convergence to these distributions. The facts and notions are used in the previous sections.
For more detailed study see \cite{GnedKolm}, \cite{IbrLin}, and \cite{Petr}.

\subsection{Self-decomposable distributions}\newparagraph
\begin{Definition}
A distribution function $F$ is called self-decomposable if for any $\alpha>1$ there exists a distribution function $F_\alpha$ such that $F(x)=\int_{\R} F_\alpha(x-y)\dd F(\alpha y)$ for all $x\in\R$.
\end{Definition}

Self-decomposable distribution functions are also called \textit{distributions functions of class} $\SelfDecLaws$ (or  $\SelfDecLaws$ \textit{functions} for short). Every self-decomposable distribution function $F$ is infinitely divisible and thus its characteristic function $f(t)=\int_{\R} e^{itx} \dd F(x)$, $t\in\R$, uniquely admits \textit{L\'evy canonical representation}:
\begin{eqnarray}\label{def_LevyRepres_Probab}
f(t)=\exp\biggl\{i\gamma t-\dfrac{\sigma^2t^2}{2}+\int\limits_{|x|>0}\biggl(e^{itx}-1-\dfrac{itx}{1+x^2}\biggr)\dd L(x) \biggr\},\quad t\in\R,
\end{eqnarray}
where $\gamma\in\R$, $\sigma^2\geqslant0$, \textit{L\'evy spectral function} $L\colon \R\setminus\{0\}\to\R$ is non-decreasing on the intervals $(-\infty,0)$ and $(0,\infty)$ and it satisfies the conditions 
\begin{eqnarray*}
&&\lim_{x\to-\infty}L(x)=\lim_{x\to\infty}L(x)=0,\\
&&\int\limits_{0<|x|<\tau}x^2 \dd L(x)<\infty \quad\text{for all}\quad\tau>0. 
\end{eqnarray*}
Only for $\SelfDecLaws$ functions its L\'evy spectral functions $L$ are continuous and have one-sided derivatives on $\R/\{0\}$, where the functions $x\mapsto x L'(x)$ are non-increasing  on the intervals $(-\infty,0)$ and $(0,\infty)$ (here $L'(x)$ is a left or right derivative of $L$ at $x$).

It is known fact that any $\SelfDecLaws$ function is unimodal (see \cite{Yam} and also \cite{Wolfe1}). Any non-degenerate $\SelfDecLaws$ function is absolutely continuous (see \cite{SatoYam}, \cite{Wolfe2} and \cite{Zolot1}).

\begin{Remark}\label{rem_StrIncSelfDec_Probab}
If $F$ is a non-degenerate $\SelfDecLaws$ function, then it strictly increases on $(\lext F, \rext F)$.
\end{Remark}
This fact follows from theorem of W. N. Hudson and H. G. Tucker (see \cite{HudTuck}), which states that a set of zeroes for density of arbitrary absolutely continuous infinitely divisible distribution function either has Lebesgue measure zero or almost surely coincides with some infinite interval.  

For example, a distribution function $D_\beta$ with the following characteristic function
\begin{eqnarray*}
	f_{D_\beta}(t)=\exp\biggl\{\beta\int\limits_0^1 \dfrac{e^{itx}-1}{x} \dd x\biggr\},\quad\beta>0,\quad t\in\R,
\end{eqnarray*}
is self-decomposable. It has canonical L\'evy representation with the triplet
\begin{eqnarray*}
	\gamma=\int\limits_0^1 \dfrac{\beta\dd x}{1+x^2}= \dfrac{\beta \pi}{4},\qquad \sigma^2=0,\qquad L(x)=\beta(\ln x) \id\bigl(x\in(0,1]\bigr).
\end{eqnarray*}
Density $\rho_\beta$ of the distribution function $D_\beta$ equals zero on $(-\infty,0]$ and it satisfies the following equation (see \cite{Hensl}, \cite{Wolfe1}, and \cite{Wolfe2})
\begin{eqnarray*}
	x \rho'_\beta(x)=(\beta-1)\rho_\beta(x)-\beta\rho_\beta (x-1), \quad x>0,
\end{eqnarray*}
with the initial condition
\begin{eqnarray*}
	\rho_\beta(x)=\dfrac{e^{-\beta\EulerCon}}{\Gamma(\beta)}\cdot x^{\beta-1}, \quad 0<x\leqslant 1,
\end{eqnarray*} 
where $\EulerCon$ is the Euler constant.

The function $\rho$, defined by the equation $\rho(x)=e^\EulerCon\rho_1(x)$, $x\in\R$, is called \textit{Dickman function}. It occupies an important place in number theory (see \cite{Tenen}). According to this, the distribution corresponding to $D_1$ is called the \textit{Dickman distribution} and the distributions corresponding to $D_\beta$, $\beta>0$, are called \textit{convolution powers of the Dickman distribution} (see \cite{Hensl}).

Self-decomposable distributions are of interest because of the following theorem (see \cite{Petr}, p. 101).
\begin{Theorem}\label{th_LimitSelfDec_Probab}
Let $(A_n)_{n\in\N}$ be a sequence and  $(B_n)_{n\in\N}$ be a positive sequence. Let $(Y_j)_{j\in\N}$ be a sequence of independent random variables satisfying the following condition
\begin{eqnarray}\label{cond_UN_Probab}
\lim_{n\to\infty}\max_{j=1,\ldots,n}\Probab(|Y_j|>x B_n)=0\quad\text{for all}\quad x>0.
\end{eqnarray}
If the ditribution functions of the following sums
\begin{eqnarray}\label{def_sums_SelfDec_Probab}
\dfrac{\sum_{j=1}^{n} Y_j -A_n}{B_n}, \quad n\in\N,
\end{eqnarray}
weakly converge to a non-degenerate distribution function $F$, then $F\in\SelfDecLaws$ and $B_n\to\infty$, $B_{n+1}/B_n\to1$ as $n\to\infty$.
\end{Theorem}

Let us formulate necessary and sufficient conditions for convergence of distribution functions of \eqref{def_sums_SelfDec_Probab} to a given $\SelfDecLaws$ function (see \cite{GnedKolm}, p. 124).

\begin{Theorem}\label{th_NessSuffCond_ConvtoSelfDec_Probab}
Let $(A_n)_{n\in\N}$ be a sequence and  $(B_n)_{n\in\N}$ be a positive sequence. Let $(Y_j)_{j\in\N}$ be a sequence of independent random variables satisfying \eqref{cond_UN_Probab}. Let $F$ be an $\SelfDecLaws$ function with triplet $(\gamma, \sigma^2, L)$ in L\'evy canonical  representation \eqref{def_LevyRepres_Probab}. For the distribution functions of the sums \eqref{def_sums_SelfDec_Probab} to converge weakly to $F$, it is necessary and sufficient to satisfy
\begin{eqnarray*}
\text{$\condAone$}&&\lim_{n\to\infty}\sum_{j=1}^{n} \Probab(Y_j>x B_n)=- L(x)\quad\text{for all}\quad x>0;\\
\text{$\condAtwo$}&&\lim_{n\to\infty}\sum_{j=1}^{n} \Probab(Y_j\leqslant x B_n)= L(x)\quad\text{for all}\quad x<0;\\
\text{$\condB$}&&\lim_{n\to\infty}\dfrac{1}{B_n} \biggl(\sum_{j=1}^{n} \Expec \Bigl[Y_j\id(|Y_j|\leqslant \tau B_n)\Bigr]-A_n\biggr)=\\ &&{}\qquad\qquad\qquad\qquad\qquad\qquad=\gamma+\int\limits_{0<|x|<\tau} \dfrac{x^3\,\dd L(x)}{1+x^2}-\int\limits_{|x|\geqslant \tau} \dfrac{x\,\dd L(x)}{1+x^2} \quad\text{for all}\quad \tau>0; \\
\text{$\condC$}&&\lim_{\tau\to0}\varliminf_{n\to\infty} \dfrac{1}{B_n^2}\sum_{j=1}^{n}\Disp \Bigl[Y_j\id(|Y_j|\leqslant \tau B_n)\Bigr]=\lim_{\tau\to0}\varlimsup_{n\to\infty} \dfrac{1}{B_n^2}\sum_{j=1}^{n}\Disp \Bigl[Y_j\id(|Y_j|\leqslant \tau B_n)\Bigr]=\sigma^2.
\end{eqnarray*}
\end{Theorem}

The convergence of \eqref{def_sums_SelfDec_Probab} still holds if we replace $(B_n)_{n\in\N}$ with any equivalent sequence. It follows from the next general lemma (see \cite{Petr} p. 21).

\begin{Lemma}\label{lm_anbnWeakConv_Probab}
Let $(c_n)_{n\in\N}$ be a sequence, $r=(r_n)_{n\in\N}$ be a positive sequence. Suppose that distribution functions  $F_n$ weakly converge to a non-degenerate distribution function $F$. Then the following assertions hold:
\begin{itemize}
\item[$(i)$] If $F_n(r_n x+c_n)$ weakly converge to non-degenerate distribution function $H$, then $H(x)=F(rx+~c)$, where $c=\lim\limits_{n\to\infty} c_n$ and $r=\lim\limits_{n\to\infty} r_n$. 
\item[$(ii)$] If $c=\lim\limits_{n\to\infty} c_n$ and $r=\lim\limits_{n\to\infty} r_n$, then $F_n(r_n x+c_n)$ weakly converge to $F(r x+c)$.
\end{itemize}
\end{Lemma}

\subsection{Stable distributions}
\begin{Definition}
A distribution function $F$ is called stable if for any $a_1>0$ and $a_2>0$ there exist $a>0$ and $b$ such that $F(a x+b)=\int_\R F(a_1 (x-y))\dd F(a_2 y)$ for all $x\in\R$.
\end{Definition}

It is well known that every stable distribution function is self-decomposable and thus infinitely divisible. The class of stable distributions consists degenerate distributions, normal distributions and non-normal $\alpha$-stable distributions, $\alpha\in(0,2)$.

The characteristic function of normal distribution function $\Phi_{\mu,\sigma^2}$ (with mean $\mu\in\R$ and variance $\sigma^2>0$) has  representation \eqref{def_LevyRepres_Probab} with the triplet $(\mu,\sigma^2,0)$, where L\'evy spectral function is identically equals zero.
Standard normal distribution function $\Phi_{0,1}$ is denoted by $\Phi$, i.e.
\begin{eqnarray}\label{conc_PhiStd_Probab}
\Phi(x)=\Phi_{\mu,\sigma^2}(\mu+\sigma x).
\end{eqnarray}

We will denote by $\StL_{\alpha, \rho, \beta, \mu}$ the stable distribution functions, which are non-degenerate and non-normal. The parametrization is choosen according to so called $A$-form of representation of their characteristic functions  (see \cite{ChrWolf} p. 10): 
\begin{eqnarray*} 
f_{\StL_{\alpha, \rho, \beta, \mu}}(t)\colonequals \exp\Bigl\{i \mu t-\rho\, \kappa_\alpha|t|^\alpha \bigl(1-i\beta\, \sign(t) \,\omega(t,\alpha)\bigr)\Bigr\},
\end{eqnarray*}
where $\alpha\in(0,2)$, $\rho>0$, $\beta\in[-1,1]$, $\mu\in\R$, and
\begin{eqnarray*} 
\kappa_\alpha\colonequals
\begin{cases}
\dfrac{\Gamma(2-\alpha)\cos(\pi\alpha/2)}{1-\alpha},&\quad \alpha\in(0, 2),\, \alpha\ne 1,\\
\pi/2,&\quad \alpha=1;
\end{cases}
\end{eqnarray*}
\begin{eqnarray*}
\omega(t,\alpha)\colonequals
\begin{cases}
\tan (\pi\alpha/2),& \quad \alpha\in(0, 2),\, \alpha\ne 1,\\
-(2/\pi)\ln|t|,& \quad \alpha= 1.
\end{cases}
\end{eqnarray*}
We will call $S_{\alpha, \beta}\colonequals S_{\alpha,1,\beta,0}$, $\alpha\in(0,2)$, $\beta\in[-1,1]$, the \textit{standard $\alpha$-stable distribution functions}, which can be obtained by the formula
\begin{eqnarray*} \label{conc_StLStd_Probab}
\StL_{\alpha, \beta}(x)=\StL_{\alpha,\rho,\beta,\mu}(\mu_1+\rho^{1/\alpha} x),
\end{eqnarray*}
where $\mu_1\colonequals \mu+ (\beta \rho \ln \rho )\,\id(\alpha=1)$.

The characteristic function of  $\StL_{\alpha, \rho, \beta, \mu}$ has  representation \eqref{def_LevyRepres_Probab} with the triplet $(\gamma,0,L)$. Here
\begin{eqnarray}\label{conc_LevySpectrFuncStable_Probab}
L(x)=\dfrac{\rho(1-\beta)}{2\,|x|^\alpha}\,\id(x<0)- \dfrac{\rho(1+\beta)}{2\, |x|^{\alpha}}\,\id(x>0), \quad x\in\R\setminus\{0\},
\end{eqnarray}
and $\gamma=\mu+\alpha \beta \rho\, I_\alpha$, where 
\begin{eqnarray*}
I_\alpha=\int\limits_{0}^{\infty} \biggr(\dfrac{1}{x^\alpha(1+x^2)}-\dfrac{\sin x}{x^2}\, \id(\alpha =1)-\dfrac{1}{x^\alpha}\,\id(1<\alpha<2)\biggl)\dd x.
\end{eqnarray*}
This integral can be explicitly computed. By formula \textbf{3.781} 1. in \cite{GradRyz}, we find
\begin{eqnarray*}
I_1=\int\limits_{0}^{\infty} \biggr(\dfrac{1}{x(1+x^2)}-\dfrac{\sin x}{x^2}\biggl)\dd x=\EulerCon-1,
\end{eqnarray*}
where $\EulerCon\approx 0,5772$ is the Euler constant. By formula \textbf{3.241} 2. in \cite{GradRyz}, we have
\begin{eqnarray*}
\int\limits_{0}^{\infty} \dfrac{x^p\dd x}{1+x^2}  = \dfrac{\pi}{2\cos (\pi p/2)}, \quad p\in(-1,1).
\end{eqnarray*}
Hence we obtain
\begin{eqnarray*}
I_\alpha=\int\limits_{0}^{\infty} \dfrac{x^{-\alpha}\dd x}{1+x^2}= \dfrac{\pi}{2 \cos(\pi \alpha/2)}\quad\text{for}\quad \alpha\in(0,1), 
\end{eqnarray*}
and also
\begin{eqnarray*}
I_\alpha=\int\limits_{0}^{\infty} \biggr(\dfrac{1}{x^\alpha(1+x^2)}-\dfrac{1}{x^\alpha}\biggl)\dd x=-\int\limits_{0}^{\infty} \dfrac{x^{2-\alpha}\dd x}{1+x^2}\dd x=\dfrac{\pi}{2 \cos(\pi \alpha/2)}\quad\text{for}\quad \alpha\in(1,2).
\end{eqnarray*}

The analytic properties of stable distributions was in detail described in the monograph  \cite{Zolot2}, the related limit theorems --- in \cite{ChrWolf}. The general reviews can be found in classic  monographs  \cite{Fell}, \cite{GnedKolm}, \cite{IbrLin}, \cite{Loeve}, and  \cite{SamTaq} but it should be take into account the paper \cite{Hall}.

The fundamental role of stable distributions is explained  by the following theorem.

\begin{Theorem}\label{th_LimitStable_Probab}
The set of distribution functions that are weak limits of the distribution functions of centered and normalized sums \eqref{def_sums_SelfDec_Probab} with independent and identically distributed random variables $Y_j$, $j\in\N$, coincides with the set of stable distribution functions.
\end{Theorem}

The criterion of convergence to standard normal distribution function has the following form (see Theorem 2.6.2 in \cite{IbrLin}, but there is a typos in the formulation: ``if'' should be changed to ``iff'').
\begin{Theorem}\label{th_ConvtoNormal_Probab}
Let $(Y_j)_{j\in\N}$ be a sequence of non-degenerate independent and identically distributed random variables. There exist a sequence $(A_n)_{n\in\N}$ and a strictly positive sequence $(B_n)_{n\in\N}$ such that the distribution functions of \eqref{def_sums_SelfDec_Probab} weakly convergence as $n\to\infty$ to $\Phi$ iff the following conditions hold:
\begin{itemize}
\item[$(i)$] $\Disp Y_1<\infty$,
\item[$(ii)$] $\Probab(|Y_1|> x)=x^{-2}\varphi(x)$,
\end{itemize}
where $\varphi$ is a SVF \footnote{slowly varing function on infinity, see Subsection \ref{subsec_RV}}. For $A_n$, $n\in\N$, it can be set $A_n=n\, \Expec Y_1$, the $B_n$, $n\in\N$, can be taken from $n\, \Disp \bigl(Y_1\id(|Y_1|<B_n)\bigr)\sim B_n^2$, $n\to\infty$.
In case $(i)$ we can set $B_n=\sqrt{n \Disp Y_1}$.
\end{Theorem}

The criterion of convergence to standard $\alpha$-stable distribution functions has the following form (see \cite{SamTaq} p. 50, \cite{Whitt} p. 114).

\begin{Theorem}\label{th_ConvtoStable_Probab}
Let $(Y_j)_{j\in\N}$ be a sequence of non-degenerate independent and identically distributed random variables. There exist a sequence $(A_n)_{n\in\N}$ and a strictly positive sequence $(B_n)_{n\in\N}$ such that the distribution functions of \eqref{def_sums_SelfDec_Probab} weakly convergence as $n\to\infty$ to $\StL_{\alpha,\beta}$, $\alpha\in(0,2)$, $\beta\in[-1,1]$,  iff both
\begin{itemize}
\item[$1)$] $\Probab(|Y_1|> x)=x^{-\alpha}\varphi(x)$,
\item[$2)$] $\lim\limits_{x\to\infty}\dfrac{\Probab(Y_1>x)}{\Probab(|Y_1|> x)}=\dfrac{1+\beta}{2}$,
\end{itemize}  
where $\varphi$ is a SVF. The $B_n$, $n\in\N$, can be taken from $\lim\limits_{n\to\infty} n\,\Probab(|Y_1|> B_n)=1$. The $A_n$, $n\in\N$, can be chosen as follows:  $A_n=0$ for $\alpha\in(1,2)$, $A_n=n\,\Expec Y_1$ for $\alpha\in(1,2)$, and 
\begin{eqnarray*}
A_n= n \Expec \bigl(Y_1 \id(Y_1\leqslant B_n)\bigr) +\beta (1-\EulerCon) B_n\quad\text{for}\quad \alpha=1.
\end{eqnarray*}
\end{Theorem}

For the case $\alpha=1$ the expression of the $A_n$, $n\in\N$, is often omitted in the literature or has often a difficult form. Let us explain the possibility of the choice of $A_n$, $n\in\N$, in Theorem \ref{th_ConvtoStable_Probab}. Indeed, by the condition $\condB$ of Theorem \ref{th_NessSuffCond_ConvtoSelfDec_Probab} for $\tau=1$ and on account of above expressions of $\gamma$ and $L$ for $\mu=0$, $\rho=1$, and $\alpha=1$, the constants $A_n$, $n\in\N$, must be satisfy
\begin{eqnarray*}
\lim_{n\to\infty}\dfrac{1}{B_n} \biggl(\sum_{j=1}^{n} \Expec Y_j\id(Y_j\leqslant B_n)-A_n\biggr)=\beta(\EulerCon-1)+\beta \biggl(\int\limits_0^1 \dfrac{x\,\dd x}{1+x^2}-\int\limits_1^\infty \dfrac{\dd x}{(1+x^2)x}\biggr).
\end{eqnarray*}
The last integrals are equal, therefore, we have the required formula for $A_n$, $n\in\N$.

\subsection{The Darling theorem}\label{subsec_Darling}\newparagraph
Suppose that common distribution function of independent and identically distributed random variables has a slowly varying summary tail, i.e. $\Probab(|Y_1|>x)=\svf(x)$, where $\svf$ is some SVF that $\svf(x)\to 0$ as $x\to\infty$. In this case the distribution functions of the sums
\begin{eqnarray*}
\dfrac{\sum_{j=1}^{n_k} Y_j -A_k}{B_k}, \quad k\in\N,
\end{eqnarray*}
can not have non-degenerate weak limit for any $(A_k)_{k\in\N}$, $(B_k)_{k\in\N}$, and $(n_k)_{k\in\N}$ (see \cite{Fell} p. 320). However, here \textit{the Darling theorem} holds (see  \cite{Darl}). We formulate it in the form, which was obtained by S. V. Nagaev and V. I. Vakhtel in the paper \cite{NagVach}. Also we restrict ourselves to the case of non-negative random variables. 

\begin{Theorem}\label{th_Darling2_Probab}
Let $(Y_j)_{j\in\N}$ be a sequence of non-degenerate independent and identically distributed non-negative random variables such that $\Probab(Y_1>x)=\svf(x)$ with some SVF $\svf$. Then  we have
\begin{eqnarray*}
\lim_{n\to\infty}\Probab\biggl(n\,\Laplsvf\Bigl(\sum_{i=1}^n Y_i\Bigr)\leqslant x\biggr)=1-e^{-x} \quad\text{for all}\quad x>0,
\end{eqnarray*}
where $\Laplsvf(y):=1-\Expec e^{-Y_1/y}$, $y>0$.
\end{Theorem}
In this theorem the function $\Laplsvf$ is continuous and strictly decreasing on $(0,\infty)$. Also, according to tauberian theorem (see \cite{Fell}, p. 447, formula (5.22)), we have $\Laplsvf(x)\sim\svf(x)$, $x\to\infty$.

\subsection{Some facts from theory of regular variation}\label{subsec_RV}\newparagraph
Let $\svf$ be a positive measurable function, defined on some $[T,\infty)$ and satisfying
$\svf(cx)/\svf(x)\to1$, $x\to\infty$ for any $c>0$. Then $\svf$ is said to be a \textit{slowly varying} at infinity (SVF for short). Here we provide some useful lemmas concerning asymptotic inversion and conjugation of such functions. 

\begin{Lemma}
Suppose that $\svf$ is a SVF. Then there exists a SVF $\svf^\#$, unique up to asymptotic equivalence, that satisfies
\begin{enumerate}
\item $\lim\limits_{x\to\infty} \svf(x)\svf^\#(x\svf(x))=1$,
\item $\lim\limits_{x\to\infty} \svf^\#(x)\svf(x\svf^\#(x))=1$,
\item $\svf^{\#\#}(x)\sim\svf(x), \quad x\to\infty$.
\end{enumerate}
\end{Lemma}

The SVF $\svf^\#$ is called \textit{de Bruijn conjugation} of $\svf$. The importance of $\svf^\#$ is explained by the following lemma (see \cite{BingGoldTeug} p. 28--29).
\begin{Lemma}\label{lm_AsymInv_deBruijn_SVF}
Let $f(x)\sim x^\tau\svf(x)^\tau$, $x\to\infty$ with $\tau>0$, $\svf$ is a SVF. There exists the function $g$ that $f(g(x))\sim g(f(x))\sim x$ and $g(x)\sim x^{1/\tau} \svf^\#(x^{1/\tau})$, $x\to\infty$, where $\svf^\#$ is the de Bruijn conjugation for $\svf$. Here $g$ is determined uniquely up to asymptotic equivalence, and one version of $g$ is $f^{-1}(x)=\inf\{y\in\R: f(y)>x\}$.
\end{Lemma}

For many cases $\svf^\#$ may be asymptotically expressed in terms of $\svf$ itself (\cite{BingGoldTeug} p. 78--79).
\begin{Lemma}\label{lm_Bekessy_SVF}
Let $\svf$ be a SVF with the de Bruijn conjugation $\svf^\#$. If $\svf(x)\sim\svf(x/\svf(x))$, $ x\to\infty$, then we have $\svf^\#(x)\sim 1/\svf(x)$, $x\to\infty$.
\end{Lemma}

Let us mention one more useful assertion (see \cite{BingGoldTeug} p. 26--27).
\begin{Lemma}\label{lm_deHaan_SVF}
Let $\svf$ be a SVF that locally integrable on $[T,\infty)$ for some $T\in\R$. Let function $\deHaansvf$ define by $\deHaansvf(x):=\int_0^x (\svf(t)/t) \dd t$. Then $\deHaansvf$ is a SVF and $\deHaansvf(x)/\svf(x)\to\infty$, $x\to\infty$.
\end{Lemma}

\section*{Acknowlegments}\newparagraph
The author wishes to thank the Chebyshev Laboratory of St. Petesburg State University for providing excellent working conditions. Also the author wishes to express his gratitude to Professor M. A. Lifshits for many helpful suggestions and comments during the preparation of the paper.

\textit{Keywords and phrases}: multivariate problems,  tensor product-type random elements, approximation, average case approximation complexity, asymptotic analysis, tractability.

\textsc{Chebyshev Laboratory, St. Petersburg State University, 14th Line, 29b, St. Petersburg, 199178 Russia}\\
\textit{E-mail address}: \texttt{alexeykhartov@gmail.com}


\begin{thebibliography}{99}
\bibitem{Adl} R. J. Adler, J. Taylor, \textit{Random Fields and Geometry}, Springer, New York, 2007.

\bibitem{BingGoldTeug} N. H. Bingham, C. M. Goldie, J. L. Teugels, \textit{Regular Variation}, Camb. Univ. Press, Cambridge, 1987.

\bibitem{ChangHa} C.-H. Chang, C.-W. Ha, \textit{The Green functions of some boundary value problems via the Bernoulli and Euler polynomials}, Arch. Math. (Basel), \textbf{76} (2001), no. 5, 360--365.

\bibitem{ChrWolf} G. Christoph, W. Wolf, \textit{Convergence Theorems with a Stable Limit Law},  Math. Research 70, Akad. Verlag, Berlin, 1992.

\bibitem{Darl}  D. A. Darling, \textit{The Influence of the maximum term in the addition of independent random variables}, Trans. Amer. Math. Soc., \textbf{73} (1952), no. 1, 95--107.

\bibitem{Fell} W. Feller, \textit{Introduction to Probability Theory and Its Applications}, vol. 2, Wiley, New York, 1971.

\bibitem{GaoHanTor}  F. Gao, J. Hanning, F. Torcaso, \textit{Integrated Brownian motions and exact $L_2$-small balls}, Ann. Probab., \textbf{31} (2003), no. 3, 1320--1337.

\bibitem{GradRyz} I. S. Gradshteyn,  I. M. Ryzhik, \textit{Tables of Integrals, Series and Products}, Elsevier, Burlington, 2007.

\bibitem{GnedKolm} B. V. Gnedenko, A. N. Kolmogorov, \textit{Limit Distributions for Sums of Independent Random Variables}, Addison-Wesley, Cambridge, 1954.

\bibitem{Hall} P. Hall, \textit{A comedy of errors: the canonical form for a stable characteristic function}, Bull. London Math. Soc., \textbf{13} (1981), no. 1, 23--27.

\bibitem{Hensl} D. Hensley, \textit{The convolution powers of the Dickman function}, J. London Math. Soc., \textbf{33} (1986), no. 3, 395--406.

\bibitem{HudTuck} W. N. Hudson, H. G. Tucker, \textit{On admissible translates of infinitely divisible distributions}, Z. Wahrsch. Verw. Gebiete, \textbf{31} (1975), 65--72. 

\bibitem{IbrLin} I. A. Ibragimov, Yu. V. Linnik, \textit{Independent and Stationary Sequences of Random Variables}, Wolters-Noordhoff, Groningen, 1971.

\bibitem{KarNazNik} A. Karol, A. Nazarov, Ya. Nikitin, \textit{Small ball probabilities for Gaussian random fields and tensor products of compact operators}, Trans. Amer. Math. Soc., \textbf{360} (2008), no. 3, 1443--1474.

\bibitem{Lifsh} M. A. Lifshits, \textit{Lectures on Gaussian Processes}, Springer, New York, 2012.

\bibitem{LifPapWoz1} M. A. Lifshits, A. Papageorgiou, H. Wo\'zniakowski, \textit{Average case tractability of non-homogeneous tensor product problems}, J. Complexity, \textbf{28} (2012), no. 5--6, 539--561.

\bibitem{LifPapWoz2} M. A. Lifshits, A. Papageorgiou, H. Wo\'zniakowski, \textit{Tractability of multi-parametric Euler and Wiener integrated processes}, Probab. Math. Stat., \textbf{32} (2012), no. 1, 131--165.

\bibitem{LifTul} M. A. Lifshits, E. V. Tulyakova, \textit{Curse of dimensionality in approximation of random fields}, Probab. Math. Stat., \textbf{26} (2006), no. 1, 97--112.

\bibitem{LifZani} M. A. Lifshits, M. Zani, \textit{Approximation complexity of additive random fields}, J. Complexity, \textbf{24} (2008), no. 3, 362--379.

\bibitem{LifZani2014} M. A. Lifshits, M. Zani, \textit{Approximation of additive random fields based on standard information: average case and probabilistic settings}, \texttt{arXiv:1402.5489}

\bibitem{LinOstr} Yu. V. Linnik, I. V. Ostrovskii, \textit{Decomposition of Random Variables and Vectors}, Transl. Math. Monog., Vol. 48, AMS, Providence, Rhode Island, 1977.

\bibitem{Loeve} M. Lo\'eve, \textit{Probability Theory I}, vol. 1, Grad. Text Math. 45, Springer-Verlag, New-York, 1977.   

\bibitem{MakGolLodPod} B. M. Makarov, M. G. Goluzina, A. A. Lodkin, A. N. Podkorytov, \textit{Selected Problems in Real Analysis}, Transl. Math. Monog., vol. 107, AMS, Providence, Rhode Island, 1992.

\bibitem{NagVach} S. V. Nagaev, V. I. Vachtel, \textit{On sums of independent random variables without power moments}, Siberian Math. J., \textbf{49} (2008), no. 6, 1091--1100.

\bibitem{NazNik} A. I. Nazarov, Ya. Yu. Nikitin,  \textit{Exact $L_2$-small ball behavior of integrated Gaussian processes and spectral asymptotics of boundary value problems}, Probab. Theory Relat. Fields, \textbf{129} (2004), no. 4, 469--494.

\bibitem{NovWoz0} E. Novak, I. H. Sloan, J. F. Traub, H. Wo\'zniakowski, \textit{Essays on the Complexity of Continuous Problems}, EMS, Z\"urich, 2009.

\bibitem{NovWoz1}  E. Novak, H. Wo\'zniakowski, \textit{Tractability of Multivariate Problems. Volume I: Linear Information}, EMS Tracts Math. 6, EMS, Z\"urich, 2008.

\bibitem{NovWoz2}  E. Novak, H. Wo\'zniakowski, \textit{Tractability of Multivariate Problems. Volume II: Standard Information for Functionals}, EMS Tracts Math. 12, EMS, Z\"urich, 2010.

\bibitem{NovWoz3}  E. Novak, H. Wo\'zniakowski, \textit{Tractability of Multivariate Problems. Volume III: Standard Information for Operators}, EMS Tracts Math. 18, EMS, Z\"urich, 2012.

\bibitem{Petr} V. V. Petrov, \textit{Limit Theorems of Probability Theory: Sequences of Independent Random Variables}, Oxford Stud. Prob. 4, Clarendon Press, Oxford, 1995.

\bibitem{Rit} K. Ritter, \textit{Average-case Analysis of Numerical Problems}, Lecture Notes in Math. No. 1733, Springer, Berlin, 2000.

\bibitem{SamTaq} G. Samorodnitsky, M. S. Taqqu, \textit{Stable Non-Gaussian Random Processes: Stochastics Models with Infinite Variance}, Chapman \& Hall, London, 1994.

\bibitem{SatoYam} K.-I. Sato, M. Yamazato, \textit{On Distribution function of class L}, Z. Wahrsch. Verw. Gebiete, \textbf{43} (1978), 273--308.


\bibitem{Tenen} G. Tenenbaum, \textit{Introduction to Analytic and Probabilistic Number Theory}, Camb. Univ. Press, Cambridge, 1995.

\bibitem{TraubWasWoz1}  J. F. Traub, G. W.  Wasilkowski, H. W\'ozniakowski, \textit{Information, Uncertainty, Complexity}, Addison-Wesley, Reading MA, 1983. 

\bibitem{TraubWasWoz2}  J. F. Traub, G. W. Wasilkowski, H. W\'ozniakowski, \textit{Information-Based Complexity}, Academic Press, New York, 1988. 

\bibitem{TraubWasWoz3}  J. F. Traub, H. W\'ozniakowski, \textit{A general theory of optimal algorithms}, Academic Press, NewYork, 1980.

\bibitem{TraubWers}  J. F. Traub, A. G. Werschulz, \textit{Complexity and Information}, Camb. Univ. Press, Cambridge, 1998.

\bibitem{Vaart} A. W. van der Vaart, \textit{Asymptotic Statistics}, Camb. Univ. Press, Cambridge, 1998.

\bibitem{WasWoz} G. W. Wasilkowski, H. W\'ozniakowski, \textit{Average case optimal algorithms in Hilbert spaces}, J. Approx. Theory, \textbf{47} (1986), 17--25.

\bibitem{Whitt} W. Whitt, \textit{Stochastic-Process Limits: An Introduction to Stochastic-Process Limits and Their Application to Queues}, Springer, New York, 2002.

\bibitem{Wolfe1} S. J. Wolfe,  \textit{On the unimodality of L functions}, Ann. Math. Statist., \textbf{42} (1971), no. 3, 912--918.

\bibitem{Wolfe2} S. J. Wolfe,   \textit{On the continuity properties of L functions}, Ann. Math. Statist., \textbf{42} (1971), no. 6, 2064--2073.

\bibitem{Yam} M. Yamazato, \textit{Unimodality of infinitely divisible distribution functions of class L}, Ann. Probab., \textbf{6} (1978), no. 4, 523--531.

\bibitem{Zolot1} V. M. Zolotarev, \textit{The analytic structure of infinitely divisible laws of class L},  Litovsk. Mat. Sb.,  \textbf{3} (1963), 123--140.

\bibitem{Zolot2} V. M. Zolotarev, \textit{One-dimensional Stable Distributions}, Transl. Math. Monog., vol. 65, AMS, Providence, Rhode Island, 1986.
\end{thebibliography}
\end{document}